\documentclass[final]{siamltex}

\usepackage{amsbsy}
\usepackage{amsmath}
\usepackage{amssymb}
\usepackage{latexsym}
\usepackage{ifthen}
\usepackage{subfigure}
\usepackage{turnstile}
\usepackage{mathtools}
\usepackage{booktabs}
\usepackage{balance}
\usepackage{paralist}
\usepackage{cite}

\newtheorem{remark}{Remark}[section]
\newtheorem{example}{Example}[section]

\def \mxc{\mathbf{x}_{\mbox{\scriptsize{c}}}}

\def \Sigmac{\Sigma_{\mbox{\scriptsize{c}}}}
\def \OverlineSigmac{\overline{\Sigma}_{\mbox{\scriptsize{c}}}}
\def \Past{\mathbb{P}_{\mathbf{\theta}^{\ast}}}
\def \East{\mathbb{E}_{\mathbf{\theta}^{\ast}}}
\def \Vap{\varepsilon}
\def \Gavg{G_{\mbox{\scriptsize{avg}}}}
\def \TGavg{\widetilde{G}_{\mbox{\scriptsize{avg}}}}
\def \Xavg{\mathbf{x}_{\mbox{\scriptsize{avg}}}}
\def \Davg{D_{\mbox{\scriptsize{avg}}}}
\def \C{\mathcal{C}}
\def \PC{\mathcal{C}^{\perp}}
\def \OL{\overline{L}}
\def \WL{\widetilde{L}}

\mathtoolsset{showonlyrefs}

\title{Distributed Linear Parameter Estimation: Asymptotically Efficient Adaptive Strategies}

\author{Soummya Kar\footnotemark[2]\ \footnotemark[4]
\and Jos\'e M.~F.~Moura\footnotemark[2]\ \footnotemark[4]
\and H. Vincent Poor\footnotemark[3] \footnotemark[5]}

\begin{document}
\maketitle

\renewcommand{\thefootnote}{\fnsymbol{footnote}}

\footnotetext[2]{Department of Electrical and Computer Engineering,
Carnegie Mellon University, Pittsburgh, PA 15213, USA (soummyak@andrew.cmu.edu, moura@ece.cmu.edu).}
\footnotetext[3]{Department of Electrical Engineering, Princeton University, Princeton, NJ 08544, USA (poor@princeton.edu).}
\footnotetext[4]{The work of S. Kar and J.M.F. Moura was partially supported by NSF grants \#~CCF-1018509 and~CCF-1011903 and by AFOSR grant \#~FA9550101291.}
\footnotetext[5]{The work of H. Vincent Poor was supported by the Office of Naval Research under grant \#~N00014-12-1-0767.}

\renewcommand{\thefootnote}{\arabic{footnote}}

\begin{abstract} The paper considers the problem of distributed adaptive linear parameter estimation in multi-agent inference networks. Local sensing model information is only partially available at the agents and inter-agent communication is assumed to be unpredictable. The paper develops a generic mixed time-scale stochastic procedure consisting of simultaneous distributed learning and estimation, in which the agents adaptively assess their relative observation quality over time and fuse the innovations accordingly. Under rather weak assumptions on the statistical model and the inter-agent communication, it is shown that, by properly tuning the consensus potential with respect to the innovation potential, the asymptotic information rate loss incurred in the learning process may be made negligible. As such, it is shown that the agent estimates are asymptotically efficient, in that their asymptotic covariance coincides with that of a centralized estimator (the inverse of the centralized Fisher information rate for Gaussian systems) with perfect global model information and having access to all observations at all times. The proof techniques are mainly based on convergence arguments for non-Markovian mixed time scale stochastic approximation procedures. Several approximation results developed in the process are of independent interest.
\end{abstract}


\section{Introduction}
\label{introduction}

\subsection{Background and Motivation}
\label{backmot} Recent advances in sensing and communication technologies have enabled the proliferation of heterogeneous sensing resources in multi-agent networks, typical examples being cyberphysical systems and distributed sensor networks. Due to the large size of these networks and the presence of geographically spread resources, distributed information processing and optimization (see, for example,~\cite{olfatisaberfaxmurray07,dimakiskarmourarabbatscaglione-11}) techniques are gaining prominence. They not only offer a robust alternative to fusion center based centralized approaches, but lead to efficient use of the network resources by distributing the computing and communication burden among the agents. A key challenge in such distributed processing involves the lack of global (sensing) model information at the local agent level. Moreover, the systems under consideration are dynamic, often leading to uncertainty in the spatial distribution of the information content. The performance of existing distributed information processing and optimization schemes (see, for example,~\cite{murray:04,Xiao05ascheme,Schizas08giannakis,KarMouraRamanan-Est,JadSanTah10,jakoveticxaviermoura-10,zamp:07,usman_tsp:07,Sayed-LMS,Giannakis-LMS,olfati:cdc09,Ram-Nedich-Siam,Lobel-Ozdaglar,Yin-1,jakoveticxaviermoura-11}) based on accurate knowledge of the sensed data statistics may suffer substantially in the face of such parametric uncertainties. This necessitates the development of adaptive schemes that learn the model parameters over time in conjunction with carrying out the desired information processing task.

Motivated by the above, in this paper we focus on the problem of distributed recursive least squares parameter estimation, in which the agents have no prior knowledge of the global sensing model and of the individual observation qualities as measured in terms of the signal to noise ratio (SNR). Our goal is to develop an adaptive distributed scheme that is asymptotically efficient, i.e., achieves the same estimation performance at each agent (in terms of asymptotic covariance) as that of a (hypothetical) centralized fusion center with perfect global model information and having access to all agents' observations at all times. To this end, we develop a \emph{consensus}+\emph{innovation} scheme, \cite{KarMouraRamanan-Est}, in which the agents collaborate by exchanging (appropriate) messages with their neighbors (consensus) and fusing the acquired information with the new local observation (innovation). Apart from the issue of optimality, the inter-agent collaboration is necessary for estimator consistency, as the local observations are generally not rich enough to guarantee global observability. Lacking prior global model and local SNR information, the innovation gains at the agents are not optimal \emph{a priori}, and the agents simultaneously engage in a distributed learning process based on past data samples with the aim of recovering the optimal gains asymptotically. Thus the distributed learning process proceeds in conjunction and interacts with the estimate update. Intuitively, the overall update scheme has the structure of a certainty-equivalent control system (see, for example,~\cite{Lai-Wei,Lai-nonlinlsad} and the references therein, in the context of parameter estimation), the key difference being the distributed nature of the learning and estimation tasks. Under rather weak assumptions on the inter-agent communication (network connectivity on \emph{average}) we show that, by properly tuning the consensus potential with respect to the innovation potential, the asymptotic information rate loss incurred in the learning process may be made negligible, and the agent estimates are asymptotically efficient in that their asymptotic covariances coincide with that of the hypothetical centralized estimator. The proper tuning of the persistent consensus and innovation potentials are necessary for this optimality, leading to a mixed time-scale stochastic procedure. In this context, we note the study of mixed time-scale stochastic procedures that arise in algorithms of the simulated annealing type (see, for example,~\cite{Gelfand-Mitter}). Apart from being distributed, our scheme technically differs from~\cite{Gelfand-Mitter} in that, whereas the additive perturbation in~\cite{Gelfand-Mitter} is a martingale difference sequence, ours is a network dependent consensus potential manifesting past dependence. In fact, intuitively, a key step in the analysis is to derive pathwise strong approximation results to characterize the rate at which the consensus term/process converges to a martingale difference process. We also emphasize that our notion of mixed time-scale is different from that of stochastic algorithms with coupling (see~\cite{Borkar-stochapp,Yin-book}), where a quickly switching parameter influences the relatively slower dynamics of another state, leading to \emph{averaged} dynamics. Mixed time scale procedures of this latter type arise in multi-scale distributed information diffusion problems; see, in particular, the paper~\cite{Krishnamurthy-Yin-consensus} that studies interactive consensus formation in Markov modulated switching networks.

We comment on the main technical ingredients of the paper. Due to the mixed time-scale behavior and the non-Markovianity (induced by the learning process that uses all past information), the stochastic procedure does not fall under the purview of standard stochastic approximation (see, for example,~\cite{Nevelson}) or distributed stochastic approximation (see, for example,~\cite{tsitsiklisbertsekasathans86,Bertsekas-survey,Kushner-dist,KarMouraRamanan-Est,Stankovic-parameter,karmoura-randomtopologynoise,Li-Feng,Huang}) procedures. As such, we develop several intermediate results on the pathwise convergence rates of mixed time-scale stochastic procedures. Some of these tools are of independent interest and general enough to be applicable to other distributed adaptive information processing problems.

We briefly summarize the organization of the rest of the paper. Section~\ref{notgraph} presents notation to be used throughout. The abstract problem formulation and the mixed time-scale distributed estimation scheme are stated and discussed in Sections~\ref{sys_model} and~\ref{alg_ADLE} respectively. The main results of the paper are stated in Section~\ref{main_res}, whereas Section~\ref{app_res} presents some intermediate convergence results on recursive stochastic schemes. The key technical ingredients concerning the asymptotic properties of the distributed learning and estimation processes are obtained in Section~\ref{sec:conv-asym}, while the main results of the paper are proved in Section~\ref{sec:proof_main_res}. Finally, Section~\ref{conclusion} concludes the paper.

\subsection{Notation}
\label{notgraph} We denote the $k$-dimensional Euclidean space by
$\mathbb{R}^{k}$. The set of reals is denoted by $\mathbb{R}$, whereas $\mathbb{R}_{+}$ denotes the non-negative reals. For $a,b\in\mathbb{R}$, we will use the notations $a\vee b$ and $a\wedge b$ to denote the maximum and minimum of $a$ and $b$ respectively. The set of $k\times k$ real matrices is denoted by $\mathbb{R}^{k\times k}$. The corresponding subspace of symmetric matrices is denoted by $\mathbb{S}^{k}$. The cone of positive semidefinite matrices is denoted by $\mathbb{S}_{+}^{k}$, whereas $\mathbb{S}_{++}^{k}$ denotes the subset of positive definite matrices. The $k\times k$ identity matrix is
denoted by $I_{k}$, while $\mathbf{1}_{k}$ and $\mathbf{0}_{k}$ denote
respectively the column vector of ones and zeros in
$\mathbb{R}^{k}$. Often the symbol $0$ is used to denote the $k\times p$ zero matrix, the dimensions being clear from the context. The operator
$\left\|\cdot\right\|$ applied to a vector denotes the standard
Euclidean $\mathcal{L}_{2}$ norm, while applied to matrices it denotes the induced
$\mathcal{L}_{2}$ norm, which is equivalent to the matrix spectral radius for symmetric
matrices. The notation $A\otimes B$ is used to denote the Kronecker product of two matrices $A$ and $B$.

Time is assumed to be discrete or slotted throughout the paper. The symbols $t$ and $s$  denote time, and $\mathbb{T}_{+}$ is the discrete index set $\{0,1,2,\cdots\}$. The parameter to be estimated belongs to a subset~$\Theta$ (generally open) of the Euclidean space $\mathbb{R}^{M}$. The true (but unknown) value of the
parameter is~$\mathbf{\theta}^{\ast}$ and a
canonical element of~$\Theta$ is~$\mathbf{\theta}$. The
estimate of~$\mathbf{\theta}^{\ast}$ at time~$t$ at agent~$n$ is
 $\mathbf{x}_{n}(t)\in\mathbb{R}^{M}$. Without
loss of generality, the initial estimate,
$\mathbf{x}_{n}(0)$, at time~$0$ at agent~$n$ is a non-random
quantity.

\textbf{Spectral graph theory}: The inter-agent communication topology may be described by an \emph{undirected} graph $G=(V,E)$, with $V=\left[1\cdots N\right]$ and~$E$  the set of agents (nodes) and communication links (edges), respectively. The unordered pair $(n,l)\in E$ if there exists an edge between nodes~$n$ and~$l$. We consider simple graphs, i.e., graphs devoid of self-loops and multiple edges. A graph is connected if there exists a path\footnote{A path between nodes $n$ and $l$ of length $m$ is a sequence
$(n=i_{0},i_{1},\cdots,i_{m}=l)$ of vertices, such that $(i_{k},i_{k+1})\in E\:\forall~0\leq k\leq m-1$.}, between each pair of nodes. The neighborhood of node~$n$ is
\begin{equation}
\label{def:omega} \Omega_{n}=\left\{l\in V\,|\,(n,l)\in
E\right\}. 
\end{equation}
Node~$n$ has degree $d_{n}=|\Omega_{n}|$ (the number of edges with~$n$ as one end point.) The structure of the graph can be described by the symmetric $N\times N$ adjacency matrix, $A=\left[A_{nl}\right]$, $A_{nl}=1$, if $(n,l)\in E$, $A_{nl}=0$, otherwise. Let the degree matrix  be the diagonal matrix $D=\mbox{diag}\left(d_{1}\cdots d_{N}\right)$. By definition, the positive semidefinite matrix $L=D-A$ is called the graph Laplacian matrix. The eigenvalues of $L$ can be ordered as $0=\lambda_{1}(L)\leq\lambda_{2}(L)\leq\cdots\leq\lambda_{N}(L)$, the eigenvector corresponding to $\lambda_{1}(L)$ being $(1/\sqrt{N})\mathbf{1}_{N}$. The multiplicity of the zero eigenvalue equals the number of connected components of the network; for a connected graph, $\lambda_{2}(L)>0$. This second eigenvalue is the algebraic connectivity or the Fiedler value of
the network; see \cite{FanChung} for detailed treatment of graphs and their spectral theory.

\section{Problem Formulation}
\label{probform}

\subsection{System Model and Preliminaries}
\label{sys_model}
Let $\mathbf{\theta}^{\ast}\in\Theta$ be an $M$-dimensional (vector) parameter that is to be estimated by a network of~$N$ agents. Throughout, we assume that all the random objects are defined on a common measurable space $\left(\Omega,\mathcal{F}\right)$ equipped with a filtration $\{\mathcal{F}_{t}\}$. For the true (but unknown) parameter value $\mathbf{\theta}^{\ast}$, probability and expectation are denoted by $\mathbb{P}_{\mathbf{\theta}^{\ast}} \left[\cdot\right]$ and $\mathbb{E}_{\mathbf{\theta}^{\ast}}\left[\cdot\right]$, respectively. All inequalities involving random variables are to be interpreted a.s.~(almost surely.)

Each agent makes i.i.d. (independent and identically distributed) observations of noisy linear functions of the parameter. The observation model for the $n$-th agent is
\begin{equation}
\label{obsmod}
\mathbf{y}_{n}(t)=H_{n}\mathbf{\theta}^{\ast}+\mathbf{\zeta}_{n}(t)
\end{equation}
where: \begin{inparaenum}[i)] \item $\left\{\mathbf{y}_{n}(t)\in\mathbb{R}^{M_{n}}\right\}$ is the observation sequence for the
$n$-th agent; and \item for each $n$, $\left\{\mathbf{\zeta}_{n}(t)\right\}$ is a zero-mean temporally i.i.d.~ noise sequence of bounded variance, such that, $\mathbf{\zeta}_{n}(t)$ is $\mathcal{F}_{t+1}$ adapted and independent of $\mathcal{F}_{t}$. Moreover, the sequences $\left\{\mathbf{\zeta}_{n}(t)\right\}$ and $\left\{\mathbf{\zeta}_{l}(t)\right\}$ are mutually uncorrelated for $n\neq l$.
\end{inparaenum}
 For most practical agent network applications, each agent observes only a subset of~$M_n$ of the components of~$\theta$, with $M_{n}\ll M$. It is then necessary for the agents to collaborate by means of occasional local inter-agent message exchanges to achieve a reasonable estimate of the parameter $\mathbf{\theta}^{\ast}$. Moreover, due to inherent uncertainties in the deployment and the sensing environment, the statistics of the observation process (i.e., of the noise) are likely to be unknown \emph{a priori}. For example, the exact observation noise variance at an agent depends on several factors beyond the control of the deployment process and should be learned over time for reasonable estimation performance. In other words, prior knowledge of the spatial distribution of the information content (e.g., which agent is more accurate than the others) may not be available, and the proposed estimation approach should be able to adaptively learn the true value of information leading to an accurate weighting of the various observation resources.

 Let $R_{n}\in\mathbb{S}_{++}^{M_{n}}$ be the true covariance of the observation noise $\mathbf{\zeta}_{n}(t)$ at agent $n$. It is well known that, given perfect knowledge of $R_{n}$ for all $n$, the best linear centralized estimator $\{\mxc(t)\}$ of $\mathbf{\theta}^{\ast}$  is asymptotically normal, i.e.,
\begin{equation}
\label{cent_est}
\sqrt{t+1}\left(\mxc(t)-\mathbf{\theta}^{\ast}\right)\Longrightarrow~\mathcal{N}\left(\mathbf{0},\Sigmac^{-1}\right),
\end{equation}
provided the matrix $\Sigmac=\sum_{n=1}^{N}H_{n}^{T}R_{n}^{-1}H_{n}$ is invertible. In case the observation process is Gaussian, the best linear estimator is optimal, and $\Sigmac$ coincides with the Fisher information rate. In general, with the knowledge of the covariance only and no other specifics about the noise distribution, the above estimate is optimal, in that no other estimate achieves smaller asymptotic covariance than $\Sigmac^{-1}$ for all distributions with covariance $R_{n}$.

The goal of this paper is to develop a distributed estimator that leads to asymptotically normal estimates with the same asymptotic covariance $\Sigmac^{-1}$ at each agent under the following constraints: (1) Each agent is aware only of its local observation model $H_{n}$ and, more importantly, (2)~the true noise covariance $R_{n}$ is not known \emph{a priori} at agent $n$ and needs to be learned from the received observation samples and exchanged messages with its neighbors over time. Recently, in~\cite{JSTSP-Kar-Moura} a distributed algorithm was introduced that leads to the desired centralized asymptotic covariance at each agent but requires full model information (i.e., all the $H_{n}$'s) and the exact covariance values $R_{n}$ at all agents. This is due to the fact that, for optimal asymptotic covariance, the approach in~\cite{JSTSP-Kar-Moura} requires an appropriate innovation gain at each agent, the latter depending on all the model matrices and noise covariances. In the absence of model and covariance information, a natural alternative is to employ a certainty-equivalence type approach in which an adaptive sequential gain refinement (learning) step is incorporated into the desired estimation task. In this paper, we show that such a learning process (see Section~\ref{alg_ADLE}) is feasible in a distributed setting and, more importantly, the coupling between the learning and parameter estimation tasks does not slow down the convergence rate (measured in terms of asymptotic covariance) of the latter to $\mathbf{\theta}^{\ast}$.

\subsection{Adaptive Distributed Estimator: Algorithm $\mathcal{ADLE}$}
\label{alg_ADLE}
The adaptive distributed linear estimator ($\mathcal{ADLE}$) involves two simultaneous update rules, namely, (1) the estimate (state) update and (2) the gain update. To formalize, let $\{\mathbf{x}_{n}(t)\}$ denote the $\{\mathcal{F}_{t}\}$ adapted sequence of estimates of $\mathbf{\theta}^{\ast}$ at agent $n$.

\textbf{Estimate Update}: The estimate update at agent $n$ then proceeds as follows:
\begin{equation}
\label{est_up}
\mathbf{x}_{n}(t+1)=\mathbf{x}_{n}(t)-\beta_{t}\sum_{l\in\Omega_{n}(t)}\left(\mathbf{x}_{n}(t)-\mathbf{x}_{l}(t)\right)+\alpha_{t}K_{n}(t)\left(\mathbf{y}_{n}(t)-H_{n}\mathbf{x}_{n}(t)\right).
\end{equation}
In the above, $\{\beta_{t}\}$ and $\{\alpha_{t}\}$ represent appropriate time-varying weighting factors for the consensus (agreement) and innovation (new observation) potentials respectively, whereas $\{K_{n}(t)\}$ is an adaptively chosen $\{\mathcal{F}_{t}\}$-adapted matrix gain process. Also, $\Omega_{n}(t)$ denotes the $\{\mathcal{F}_{t+1}\}$-adapted time-varying random neighborhood of agent $n$ at time $t$.

\textbf{Gain Update}: The adaptive gain update at sensor $n$ involves another $\{\mathcal{F}_{t}\}$ adapted distributed learning process that proceeds in parallel with the estimate update. In particular, we set
\begin{equation}
\label{gain_up}
K_{n}(t)=\left(G_{n}(t)+\gamma_{t}I_{M}\right)^{-1}H_{n}^{T}\left(Q_{n}(t)+\gamma_{t}I_{M_{n}}\right)^{-1}
\end{equation}
where $\{\gamma_{t}\}$ is a sequence of positive reals, such that $\gamma_{t}\rightarrow 0$ as $t\rightarrow\infty$, and the positive semidefinite matrix sequences $\{Q_{n}(t)\}$ and $\{G_{n}(t)\}$ evolve as follows:
\begin{equation}
\label{Q_up}
Q_{n}(t+1)=\frac{1}{t}\sum_{s=0}^{t}\mathbf{y}_{n}(s)\mathbf{y}_{n}^{T}(s)-\left(\frac{1}{t}\sum_{s=0}^{t-1}\mathbf{y}_{n}(s)\right)\left(\frac{1}{t}\sum_{s=0}^{t-1}\mathbf{y}_{n}(s)\right)^{T},
\end{equation}
and
\begin{equation}
\label{G_up}
G_{n}(t+1)=G_{n}(t)-\beta_{t}\sum_{l\in\Omega_{n}(t)}\left(G_{n}(t)-G_{l}(t)\right)+\alpha_{t}\left(H_{n}^{T}\left(Q_{n}(t)+\gamma_{t}I_{M_{n}}\right)^{-1}H_{n}-G_{n}(t)\right)
\end{equation}
with positive semidefinite initial conditions $Q_{n}(0)$ and $G_{n}(0)$ respectively.

Before discussing further, we formalize assumptions on the model, the time-varying communication topology and the algorithm weight sequences $\{\alpha_{t}\}$ and $\{\beta_{t}\}$ in the following:

\textbf{(A.1)}: \emph{The observation model is globally observable, i.e., the (normalized) Grammian matrix
\begin{equation}
\label{glob_obs}\OverlineSigmac=\frac{1}{N}\sum_{n=1}^{N}H_{n}^{T}R_{n}^{-1}H_{n}
\end{equation}
is invertible, where $R_{n}$ denotes the non-singular true (but unknown) covariance of the observation noise $\mathbf{\zeta}_{n}(t)$ at agent $n$.}

\textbf{(A.2)}: \emph{The $\{\mathcal{F}_{t+1}\}$-adapted sequence $\{L_{t}\}$ of communication network Laplacians (modeling the agent communication neighborhoods $\left\{\Omega_{n}(t)\right\}$ at each time $t$) is temporally i.i.d. with $L_{t}$ being independent of $\mathcal{F}_{t}$ for each $t$. Further, the sequence $\{L_{t}\}$ is connected on the average, i.e., $\lambda_{2}(\overline{L})>0$, where $\overline{L}=\East[L_{t}]$ denotes the mean Laplacian.}

\textbf{(A.3)}: \emph{The sequences $\left\{L_{t}\right\}$ and $\left\{\mathbf{\zeta}_{n}(t)\right\}_{n\in V}$ are mutually independent.}

\textbf{(A.4)}: \emph{There exists $\varepsilon_{1}>0$, such that for all $n$, $\East\left[\|\mathbf{\zeta}_{n}(t)\|^{2+\varepsilon_{1}}\right]<\infty$.}

\textbf{(A.5)}: \emph{The weight sequences $\{\alpha_{t}\}$ and $\{\beta_{t}\}$ are given by
\begin{equation}
\label{weight}
\alpha_{t}=\frac{a}{(t+1)^{\tau_{1}}}~~~\mbox{and}~~~\beta_{t}=\frac{b}{(t+1)^{\tau_{2}}},
\end{equation}
where $a, b>0$, $0<\tau_{2}\leq \tau_{1}\leq 1$ and $\tau_{1}>\tau_{2}+1/(2+\varepsilon_{1})+1/2$.}
\begin{remark}
\label{rem:ass} \textup{Note that the global observability requirement in~\textbf{(A.1)} is quite weak and, in fact, is necessary to attain estimator consistency in a centralized setting. In a sense, the assumption~\textbf{(A.1)} on the global sensing model and the \emph{mean connectivity} condition in~\textbf{(A.2)} provide minimal structural conditions for attaining \emph{distributed observability}, i.e., the ability to obtain consistent parameter estimates in the proposed distributed information setting. Intuitively, the necessity of \textbf{(A.2)} (in addition to~\textbf{(A.1)}) for such distributed observability stems from the observation that, in general, in the absence of local observability a disconnected inter-agent communication network would lead to multiple communication-disjoint agent components, none with sufficient informative measurements to consistently estimate the true parameter. We emphasize that the \emph{mean} network connectivity assumption formalized in~\textbf{(A.2)}, which generalizes the notion of connectivity in static communication topologies to dynamic stochastic scenarios, models a broad class of agent networks with unpredictable communication; for instance,~\textbf{(A.2)} allows for spatially correlated communication link failures (often resulting from multi-agent interference) and subsumes the commonly used packet erasure model in gossip type of agent communications~\cite{Boyd-GossipInfTheory}. On the other hand, in the current setting, we assume that the inter-agent communication is noise-free and unquantized in the event of an active communication link; the problem of quantized data exchange in networked control systems (see, for example,~\cite{tm04,ms03,Li-Baillieul,KarMouraRamanan-Est}) is an active research topic.}

\textup{We comment on the choice of the weight sequences $\{\beta_{t}\}$ and $\{\alpha_{t}\}$ associated with the consensus and innovation potentials respectively. From~\textbf{(A.5)} we note that both the excitations for agent-collaboration (consensus) and local innovation are persistent, i.e., the sequences $\{\beta_{t}\}$ and $\{\alpha_{t}\}$ sum to $\infty$ -  a standard requirement in stochastic approximation type algorithms to drive the updates to the desired limit from arbitrary initial conditions. Further, the square summability of $\{\alpha_{t}\}$ ($\tau_{1}>1/2$) is required to mitigate the effect of stochastic sensing noise perturbing the innovations. The requirement $\beta_{t}/\alpha_{t}\rightarrow\infty$ as $t\rightarrow\infty$ ($\tau_{1}>\tau_{2}$), i.e., the asymptotic domination of the consensus potential over the local innovations ensures the right information mixing thus, as shown below, leading to optimal estimation performance. Technically, the different asymptotic decay rates of the two potentials lead to mixed time-scale stochastic recursions whose analyses require new techniques in stochastic approximation as developed in the paper.}
\end{remark}

\begin{figure}
         \centering
        \includegraphics[width=.5\linewidth]{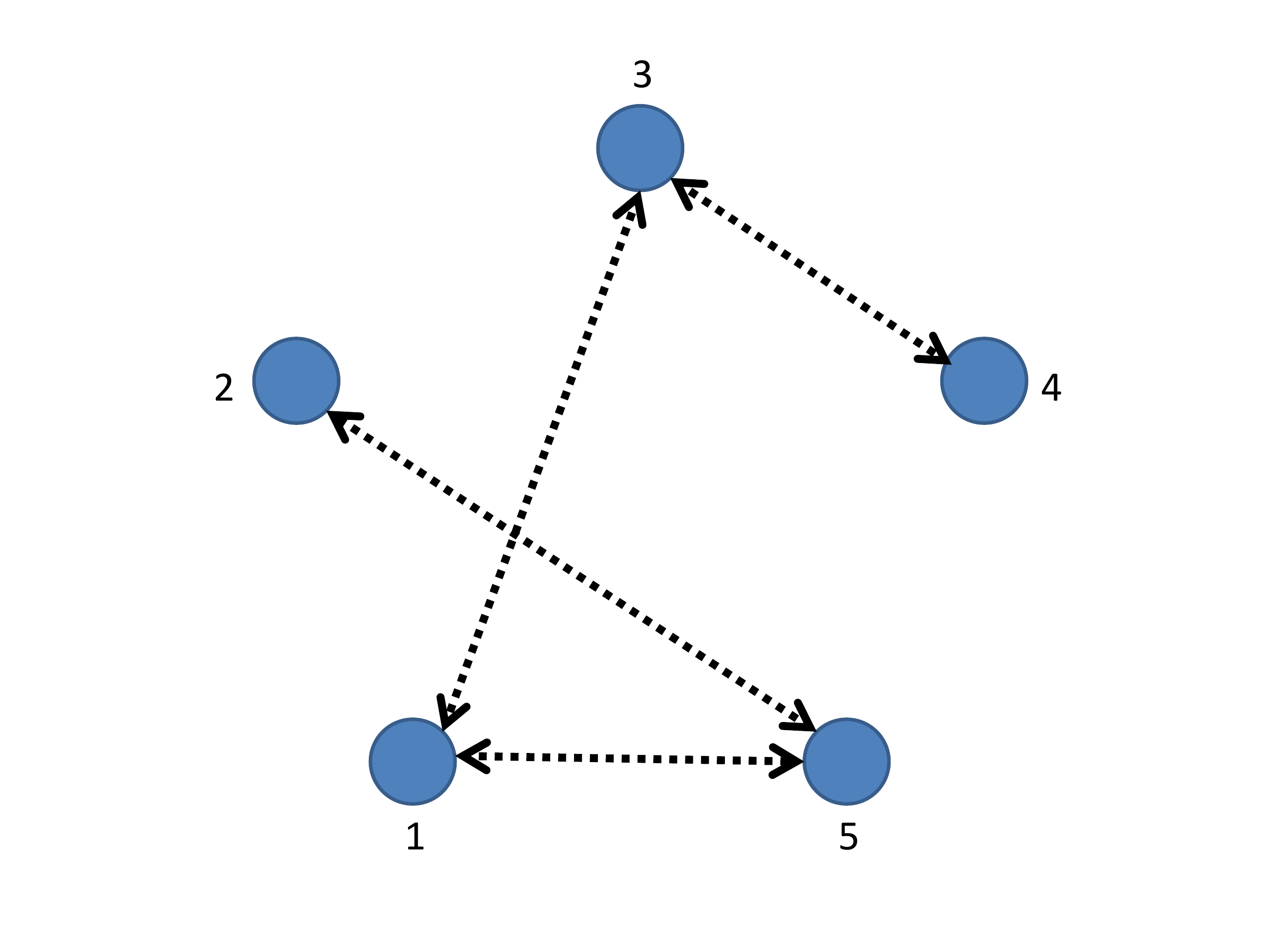}
        \caption{An example: circles depict agents and dotted lines bi-directional communication links.}
        \label{fig:ex}
\end{figure}

\begin{example}
\label{ex:alg} \textup{As an illustration, consider the agent model in Fig.~\ref{fig:ex} with $N=5$ agents. The vector parameter $\mathbf{\theta}^{\ast}\in\mathbb{R}^{5}$ in this example may have a physical interpretation, for example, with $\theta^{\ast}_{n}$, the $n$-th component of $\mathbf{\theta}^{\ast}$, indicating the (unknown) intensity of a source geographically co-located with agent $n$, $n=1,\cdots,5$. Each agent $n$ observes a scalar sequence
\begin{equation}
\label{ex:alg1}y_{n}(t)=\left(\theta^{\ast}_{n-1}+\theta^{\ast}_{n}+\theta^{\ast}_{n+1}\right)+\zeta_{n}(t),
\end{equation}
perhaps corresponding to a superposition of local source intensities, where we adopt the convention that $\theta^{\ast}_{0}=\theta^{\ast}_{5}$ and $\theta^{\ast}_{6}=\theta^{\ast}_{1}$. It is readily seen that the local agent observations are not globally observable for $\mathbf{\theta}^{\ast}$. In fact, in this example, without collaboration, no agent $n$ would be able to reconstruct even the local intensity $\theta_{n}$. The collective observation model is however globally observable for $\mathbf{\theta}^{\ast}$, i.e.,~\textbf{(A.1)} holds. The dotted lines denote the potential inter-agent communication links (possibly switching stochastically between on and off) through which the locally unobservable agents may collaborate by information exchange. By abstracting the above model in terms of the generic notation, the $\mathcal{ADLE}$ estimate update rule~\eqref{est_up} at an agent $n$, say $n=3$, then takes the form
\begin{align}
\label{ex:alg2}\mathbf{x}_{3}(t+1){}&=\mathbf{x}_{3}(t)-\beta_{t}\left(2\mathbf{x}_{3}(t)-\mathbf{x}_{1}(t)-\mathbf{x}_{4}(t)\right)\\
{}&
+\alpha_{t}\left(G_{3}(t)+\gamma_{t}I_{5}\right)^{-1}H_{3}^{T}\left(Q_{3}(t)+\gamma_{t}\right)^{-1}\left(y_{3}(t)-x_{3,2}(t)-x_{3,3}(t)-x_{3,4}(t)\right),
\end{align}
where $H_{3}=[0~1~1~1~0]$, $Q_{3}(t)$ denotes the (scalar) sample covariance~\eqref{Q_up}, $x_{3,l}(t)$ denotes the $l$-th component of $\mathbf{x}_{3}(t)$ with $l=2,3,4$, and $G_{3}(t)$ is updated as
\begin{equation}
\label{ex:alg3}
G_{3}(t+1)=G_{3}(t)-\beta_{t}\left(2G_{3}(t)-G_{1}(t)-G_{4}(t)\right)+\alpha_{t}\left(H_{3}^{T}\left(Q_{3}(t)+\gamma_{t}\right)^{-1}H_{3}-G_{3}(t)\right).
\end{equation}
In the above, we assumed that at time $t$, both the communication links $(1,3)$ and $(3,4)$ are active. Assuming that the stochasticity, if any, in the link formation satisfies~\textbf{(A.2)}, the following analysis will show that the above estimate sequences will \emph{optimally} converge to $\mathbf{\theta}^{\ast}$ a.s. as $t\rightarrow\infty$.}
\end{example}

\section{Main Results}
\label{main_res} We formally state the main results of the paper, the proofs being provided in Section~\ref{sec:proof_main_res}.

The first result concerns the asymptotic agreement or consensus among the various agent estimates.
\begin{theorem}
\label{th:estc} Let assumptions~\textbf{(A.1)}-\textbf{(A.5)} hold. Then for each $\tau_{0}$ such that
\begin{equation}
\label{th:estc1}
0\leq\tau_{0}<\tau_{1}-\tau_{2}-\frac{1}{2+\Vap_{1}},
\end{equation}
we have
\begin{equation}
\label{th:estc2}
\Past\left(\lim_{t\rightarrow\infty}(t+1)^{\tau_{0}}\left\|\mathbf{x}_{n}(t)-\mathbf{x}_{l}(t)\right\|=0\right)=1
\end{equation}
for any pair of agents $n$ and $l$.
\end{theorem}

Theorem~\ref{th:estc} relates the rate of agreement to the difference $\tau_{1}-\tau_{2}$ of the algorithm weight parameters, the latter quantifying the relative intensities of the \emph{global agreement} and \emph{local innovation} potentials. Notably, the order of this convergence is independent of the network topology (as long as it is connected in the mean) and the distributed gain learning process~\eqref{gain_up}-\eqref{G_up}. In fact, as will be evident from the proof arguments, the local covariance learning step in~\eqref{Q_up} may be replaced by any other consistent learning procedure, still retaining the order of convergence in Theorem~\ref{th:estc}.
\begin{theorem}
\label{th:estcons} Let assumptions~\textbf{(A.1)}-\textbf{(A.5)} hold with $\tau_{1}=1$ and $a\geq 1$. Then, for each $n$ the estimate sequence $\{\mathbf{x}_{n}(t)\}$ is strongly consistent. In particular, we have
\begin{equation}
\label{th:estcons1}
\Past\left(\lim_{t\rightarrow\infty}(t+1)^{\tau}\left\|\mathbf{x}_{n}(t)-\mathbf{\theta}^{\ast}\right\|=0\right)=1
\end{equation}
for each $n$ and $\tau\in [0,1/2)$.
\end{theorem}

The consistency in Theorem~\ref{th:estcons} is order optimal in that~\eqref{th:estcons1} fails to hold (unless the noise covariances are degenerate) with an exponent $\tau\geq 1/2$ for any (including centralized) estimation procedure, which is due to the fact that the optimal (centralized) estimator is asymptotically normal with non-degenerate (positive definite) asymptotic covariance.


The next result concerns the asymptotic efficiency of the estimates generated by the distributed $\mathcal{ADLE}$.
\begin{theorem}
\label{th:estn}Let assumptions~\textbf{(A.1)}-\textbf{(A.5)} hold with $\tau_{1}=1$ and $a=1$, and let $\Sigmac=N\OverlineSigmac$. Then, for each $n$
\begin{equation}
\label{th:estn200}
\sqrt(t+1)\left(\mathbf{x}_{n}(t)-\mathbf{\theta}^{\ast}\right)\Longrightarrow\mathcal{N}\left(\mathbf{0},\Sigmac^{-1}\right),
\end{equation}
where $\mathcal{N}(\cdot,\cdot)$ and $\Longrightarrow$ denote the Gaussian distribution and weak convergence, respectively.
\end{theorem}

Referring to the introductory discussion in Section~\ref{sys_model}, we note that the distributed and adaptive $\mathcal{ADLE}$ achieves the best error covariance decay among the class of linear centralized estimators and is optimal in the Fisher information sense if the noise process is Gaussian. In a sense, Theorem~\ref{th:estn} reinforces the applicability
and advantage of distributed estimation schemes. Apart from issues of robustness, implementing a
centralized estimator is communication-intensive as it requires transmitting all sensor data to
a fusion center at all times. On the other hand, the distributed $\mathcal{ADLE}$ algorithm involves only sparse local
communication among the sensors at each step, and achieves the performance of a centralized estimator
asymptotically as long as the communication network stays connected in the mean. Further, note that the assumption $a=1$ is not necessary for asymptotic normality of the $\mathcal{ADLE}$ estimates; however, the optimality (asymptotic efficiency) is no longer guaranteed for $a\neq 1$, i.e., the resulting asymptotic covariance of the estimates deviate from $\Sigmac^{-1}$.

\section{Some Approximation Results}
\label{app_res}
In this section we establish several strong (pathwise) convergence results for generic mixed time-scale stochastic recursive procedures (the proofs being provided in Appendix~\ref{app:app_res}). These are of independent interest and will be used in subsequent sections to analyze the properties of the $\mathcal{ADLE}$ scheme.

Throughout this section, by $\{\mathbf{z}_{t}\}$, we will denote an $\{\mathcal{F}_{t}\}$ adapted stochastic process taking values in some Euclidean space or some subset of symmetric matrices. The initial condition $\mathbf{z}_{0}$ will be assumed to be deterministic unless otherwise stated. Further, the probability space is assumed to be rich enough to allow the definition of various auxiliary processes governing the recursive evolution of $\{\mathbf{z}_{t}\}$. Since the results in this section concern generic stochastic processes not necessarily tied to the parameter vector, the $\mathbf{\theta}^{\ast}$ indexing in the probability and expectation will be dropped temporarily.

We start by quoting a convergence rate result from~\cite{JSTSP-Kar-Moura} on deterministic recursions with time-varying coefficients.
\begin{lemma}[Lemmas~4 and~5 of~\cite{JSTSP-Kar-Moura}]
\label{lm:JSTSP-det}
Let $\{\mathbf{z}_{t}\}$ be an $\mathbb{R}_{+}$ valued sequence satisfying
\begin{equation}
\label{lm:JSTSP-det1}
\mathbf{z}_{t+1}\leq (1-r_{1}(t))\mathbf{z}_{t}+r_{2}(t),
\end{equation}
where $\{r_{1}(t)\}$ and $\{r_{2}(t)\}$ are deterministic sequences with
\begin{equation}
\label{lm:JSTSP-det2}
\frac{a_{1}}{(t+1)^{\delta_{1}}}\leq r_{1}(t)\leq 1~~\mbox{and}~~r_{2}(t)\leq \frac{a_{2}}{(t+1)^{\delta_{2}}},
\end{equation}
and $a_{1}>0$, $a_{2}>0$, $0\leq \delta_{1}\leq 1$, $\delta_{2}>0$. Then, if $\delta_{1}<\delta_{2}$, $(t+1)^{\delta_{0}}\mathbf{z}_{t}\rightarrow 0$ as $t\rightarrow\infty$, for all $0\leq\delta_{0}<\delta_{2}-\delta_{1}$. Also, if $\delta_{1}=\delta_{2}$, the sequence $\{\mathbf{z}_{t}\}$ remains bounded, i.e.,
$\sup_{t\geq 0}\|\mathbf{z}_{t}\|<\infty$.
\end{lemma}

We now develop a stochastic analog of Lemma~\ref{lm:JSTSP-det} in which the weight sequence $\{r_{1}(t)\}$ is a random process with some mixing conditions.
\begin{lemma}
\label{lm:seq-gen}
Let $\{\mathbf{z}_{t}\}$ be an $\{\mathcal{F}_{t}\}$ adapted $\mathbb{R}_{+}$ valued process satisfying
\begin{equation}
\label{lm:JSTSP1}
\mathbf{z}_{t+1}\leq (1-r_{1}(t))\mathbf{z}_{t}+r_{2}(t).
\end{equation}
In the above, $\{r_{1}(t)\}$ is an $\{\mathcal{F}_{t+1}\}$ adapted process, such that for all $t$, $r_{1}(t)$ satisfies $0\leq r_{1}(t)\leq 1$ and
\begin{equation}
\label{lm:JSTSP2}
\frac{a_{1}}{(t+1)^{\delta_{1}}}\leq\mathbb{E}\left[r_{1}(t)~|~\mathcal{F}_{t}\right]\leq 1
\end{equation}
with $a_{1}>0$ and $0\leq \delta_{1}\leq 1$. The sequence $\{r_{2}(t)\}$ is deterministic, $\mathbb{R}_{+}$ valued, and satisfies $r_{2}(t)\leq a_{2}/(t+1)^{\delta_{2}}$ with $a_{2}>0$ and $\delta_{2}>0$. Then, if $\delta_{1}<\delta_{2}$,
$(t+1)^{\delta_{0}}\mathbf{z}_{t}\rightarrow 0$ as $t\rightarrow\infty$ for all $0\leq\delta_{0}<\delta_{2}-\delta_{1}$.
\end{lemma}

Versions of Lemma~\ref{lm:seq-gen} with stronger assumptions on the weight sequences were used in earlier work. For example, the deterministic version (Lemma~\ref{lm:JSTSP-det}) was proved in~\cite{KarMouraRamanan-Est}, whereas a version with i.i.d.~weight sequences was used in~\cite{JSTSP-Kar-Moura}. Further, several variants under varying assumptions exist in the literature based on generalized supermartingale convergence theorems; see for example~\cite{Polyak-book,tsitsiklisbertsekasathans86,Bertsekas-survey}. However, for reasons to be clear soon, in this work there will be instances in which the memoryless assumption on the weight sequences is too restrictive. Hence, we develop the version stated in Lemma~\ref{lm:seq-gen}.

The following result will be used to quantify the rate of convergence of distributed vector or matrix valued recursions to their network-averaged behavior.
\begin{lemma}
\label{lm:mean-conv} Let $\{\mathbf{z}_{t}\}$ be an $\mathbb{R}_{+}$ valued $\{\mathcal{F}_{t}\}$ adapted process that satisfies
\begin{equation}
\label{lm:mean-conv1}
\mathbf{z}_{t+1}\leq \left(1-r_{1}(t)\right)\mathbf{z}_{t}+r_{2}(t)U_{t}\left(1+J_{t}\right).
\end{equation}
Let the weight sequences $\{r_{1}(t)\}$ and $\{r_{2}(t)\}$ satisfy the hypothesis of Lemma~\ref{lm:seq-gen}.
Further, let $\{U_{t}\}$ and $\{J_{t}\}$ be $\mathbb{R}_{+}$ valued $\{\mathcal{F}_{t}\}$ and $\{\mathcal{F}_{t+1}\}$ adapted processes respectively with $\sup_{t\geq 0}\|U_{t}\|<\infty$ a.s. The process $\{J_{t}\}$ is i.i.d.~with $J_{t}$ independent of $\mathcal{F}_{t}$ for each $t$ and satisfies the moment condition $\mathbb{E}\left[\left\|J_{t}\right\|^{2+\varepsilon_{1}}\right]<\kappa<\infty$ for some $\varepsilon_{1}>0$ and a constant $\kappa>0$. Then, for every $\delta_{0}$ such that
\begin{equation}
\label{lm:mean-conv5}
0\leq\delta_{0}<\delta_{2}-\delta_{1}-\frac{1}{2+\varepsilon_{1}},
\end{equation}
we have $(t+1)^{\delta_{0}}\mathbf{z}_{t}\rightarrow 0$ a.s. as $t\rightarrow\infty$.
\end{lemma}

The key difference between Lemma~\ref{lm:mean-conv} and Lemma~\ref{lm:seq-gen} is that the processes associated with the sequence $\{r_{2}(t)\}$ are now stochastic.

\begin{lemma}
\label{lm:conn} Let $\{\mathbf{z}_{t}\}$ be an $\mathbb{R}^{NM}$ valued $\{\mathcal{F}_{t}\}$ adapted process such that $\mathbf{z}_{t}\in\mathcal{C}^{\perp}$ (see~\eqref{prop:est3} in Appendix~\ref{app:est} for the definition of the consensus subspace $\mathcal{C}$ and its orthogonal complement $\mathcal{C}^{\perp}$) for all $t$. Also, let $\{L_{t}\}$ be an i.i.d.~sequence of Laplacian matrices as in assumption~\textbf{(A.2)} that satisfies
\begin{equation}
\label{Lap_cond}
\lambda_{2}(\overline{L})=\lambda_{2}\left(\mathbb{E}[L_{t}]\right)>0,
\end{equation}
with $L_{t}$ being $\mathcal{F}_{t+1}$ adapted and independent of $\mathcal{F}_{t}$ for all $t$. Then there exists a measurable $\{\mathcal{F}_{t+1}\}$ adapted $\mathbb{R}_{+}$ valued process $\{r_{t}\}$ (depending on $\{\mathbf{z}_{t}\}$ and $\{L_{t}\}$) and a constant $c_{r}>0$, such that $0\leq r_{t}\leq 1$ a.s. and
\begin{equation}
\label{lm:conn20}
\left\|\left(I_{NM}-\beta_{t}L_{t}\otimes I_{M}\right)\mathbf{z}_{t}\right\|\leq\left(1-r_{t}\right)\left\|\mathbf{z}_{t}\right\|
\end{equation}
with
\begin{equation}
\label{lm:conn2}
\mathbb{E}\left[r_{t}~|~\mathcal{F}_{t}\right]\geq\frac{c_{r}}{(t+1)^{\tau_{2}}}~~\mbox{a.s.}
\end{equation}
for all $t$ large enough, where the weight sequence $\{\beta_{t}\}$ and $\tau_{2}$ are defined in~\eqref{weight}.
\end{lemma}
\begin{remark}
\label{rem:Lap} We comment on the necessity of the various technicalities  in the statement of Lemma~\ref{lm:conn}. Let $\mathcal{P}_{NM}$ denote the matrix $(1/N)\left(\mathbf{1}_{N}\otimes I_{M}\right)\left(\mathbf{1}_{N}\otimes I_{M}\right)^{T}$ and $\mathcal{P}_{NM}\mathbf{z}_{t}=\mathbf{0}$ since $\mathbf{z}_{t}\in\mathcal{C}^{\perp}$. With this, a naive approach of showing the existence of such a process $\{r_{t}\}$ would be to use the submultiplicative inequality
\begin{equation}
\label{lm:conn3}
\left\|\left(I_{NM}-\beta_{t}L_{t}\otimes I_{M}-\mathcal{P}_{NM}\right)\mathbf{z}_{t}\right\|\leq\left\|\left(I_{NM}-\beta_{t}L_{t}\otimes I_{M}-\mathcal{P}_{NM}\right)\right\|\left\|\mathbf{z}_{t}\right\|.
\end{equation}
Using properties of the Laplacian and the matrix $\mathcal{P}_{NM}$, it can be shown that for sufficiently large $t$
\begin{equation}
\label{lm:conn4}
\left\|\left(I_{NM}-\beta_{t}L_{t}\otimes I_{M}-\mathcal{P}_{NM}\right)\mathbf{z}_{t}\right\|\leq\left(1-\beta_{t}\lambda_{2}(L_{t})\right)\left\|\mathbf{z}_{t}\right\|.
\end{equation}
With this we may choose to define the desired sequence $\{r_{t}\}$ in Lemma~\ref{lm:conn} by
\begin{equation}
\label{lm:conn5}
r_{t}=\beta_{t}\lambda_{2}(L_{t})
\end{equation}
for all $t$. Indeed, $\{r_{t}\}$ thus defined satisfies $0\leq r_{t}\leq 1$ and~\eqref{lm:conn} (at least for $t$ large enough). Since $L_{t}$ is independent of $\mathcal{F}_{t}$, we obtain
\begin{equation}
\label{lm:conn6}
\mathbb{E}[\lambda_{2}(L_{t})~|~\mathcal{F}_{t}]=\mathbb{E}[\lambda_{2}(L_{t})]\leq\lambda_{2}(\overline{L}),
\end{equation}
where the last inequality is a consequence of Jensen's inequality applied to the concave functional $\lambda_{2}(\cdot)$. Thus the hypothesis $\lambda_{2}(\overline{L)}>0$ does not shed any light to whether $\mathbb{E}[\lambda_{2}(L_{t})]>0$ or not. Unfortunately, it turns out that in the gossip type of communication setting, in which none of the network instances are connected, $\lambda_{2}(L_{t})=0$ a.s. Hence, in such cases $\mathbb{E}[\lambda_{2}(L_{t})]$ is actually 0. This in turn implies that the $\{r_{t}\}$ proposed in~\eqref{lm:conn5} violates the requirement~\eqref{lm:conn2} of Lemma~\ref{lm:conn}. This necessitates an altogether different approach for constructing the desired sequence $\{r_{t}\}$. As shown in the following, such an $r_{t}$ is no longer independent of $\mathcal{F}_{t}$, being a function of both $L_{t}$ and $\mathbf{z}_{t}$ in general.
\end{remark}

\section{Convergence and Asymptotic Properties}
\label{sec:conv-asym} In this section we establish asymptotic properties of the $\mathcal{ADLE}$ and the associated distributed learning and estimation processes. The key technical result concerning the adaptive gain learning process is presented in Lemma~\ref{lm:gain_conv}, whereas, the major convergence properties of the estimate processes concerning boundedness of the agent estimates, inter-agent estimate consensus and estimate consistency are obtained in Lemma~\ref{lm:estb}, Lemma~\ref{lm:estc} and Lemma~\ref{lm:avgr} respectively. The assumptions~\textbf{(A.1)}-\textbf{(A.5)} are assumed to hold throughout the section.

The following result concerns the convergence of the online gain approximation processes $\{K_{n}(t)\}$ to their optimal counterparts $K_{n}=\OverlineSigmac^{-1}H_{n}^{T}R_{n}^{-1}$.
\begin{lemma}
\label{lm:gain_conv} For each $n$ the gain sequence $\{K_{n}(t)\}$ (given by~\eqref{gain_up}-\eqref{G_up}) converges to $K_{n}=\OverlineSigmac^{-1}H_{n}^{T}R_{n}^{-1}$ a.s., i.e.,
\begin{equation}
\label{lm:gain_conv1}\Past\left(\lim_{t\rightarrow\infty}K_{n}(t)=\OverlineSigmac^{-1}H_{n}^{T}R_{n}^{-1}\right)=1.
\end{equation}
\end{lemma}

The proof is accomplished in terms of several intermediate steps that highlight the interaction between the dynamics of distributed collaboration and local adaptation. To this end, we first investigate the processes $\{G_{n}(t)\}$; see~\eqref{G_up}. The processes $\{G_{n}(t)\}$ may be viewed as approximations of the normalized Grammian and, as will be shown in the following, converge to $\OverlineSigmac$. The following assertion concerns the consensus of the approximate Grammians to their network average and is stated as follows (see Appendix~\ref{app:est} for a proof):
\begin{lemma}
\label{lm:Gcons}
For each $n$ we have
\begin{equation}
\label{lm:Gcons1}
\Past\left(\lim_{t\rightarrow\infty}\left\|G_{n}(t)-\Gavg(t)\right\|=0\right)=1,
\end{equation}
where $\Gavg(t)=\frac{1}{N}\sum_{n=1}^{N}G_{n}(t)$ is the instantaneous network-averaged Grammian.
\end{lemma}
On the basis of Lemma~\ref{lm:Gcons}, to show the convergence of the approximate (normalized) Grammian sequences to $\OverlineSigmac$, it suffices to show the convergence of the network-averaged sequence $\{\Gavg(t)\}$ to the latter. This is obtained in the following lemma (see Appendix~\ref{app:est} for a proof).
\begin{lemma}
\label{lm:Gavg}
The following holds:
\begin{equation}
\label{lm:Gavg1}
\Past\left(\lim_{t\rightarrow\infty}\Gavg(t)=\OverlineSigmac\right)=1.
\end{equation}
\end{lemma}
We now complete the proof of Lemma~\ref{lm:gain_conv}.
\begin{proof}[Proof of Lemma~\ref{lm:gain_conv}] It follows from Lemma~\ref{lm:Gcons} and Lemma~\ref{lm:Gavg} that
\begin{equation}
\label{lm:gain_conv4}
\Past\left(\lim_{t\rightarrow\infty}G_{n}(t)=\OverlineSigmac\right)=1
\end{equation}
for all $n=1,\cdots,N$. The assertion in Lemma~\ref{lm:gain_conv} is immediate from~\eqref{lm:gain_conv4} and the observation that $Q_{n}(t)\rightarrow R_{n}$ and $\gamma_{t}\rightarrow 0$ as $t\rightarrow\infty$.
\end{proof}

The rest of the section is concerned with the convergence analysis of the estimate sequences $\{\mathbf{x}_{n}(t)\}$ generated by the $\mathcal{ADLE}$. Several results on the convergence behavior of the estimates are presented, culminating in the proofs of the main results of the paper in Section~\ref{sec:proof_main_res}.
\begin{lemma}
\label{lm:estb}The estimate sequences $\{\mathbf{x}_{n}(t)\}$ generated by the $\mathcal{ADLE}$ algorithm (see~\eqref{est_up}) are pathwise bounded, i.e., for each $n$, $\sup_{t\geq 0}\left\|\mathbf{x}_{n}(t)\right\|<\infty$ a.s.
\end{lemma}

The proof of this Lemma involves a Lyapunov type argument. The decay rate estimates obtained in the next two Propositions (see Appendix~\ref{app:est} for proofs) are associated with certain time-varying spectral operators. They will be used in the construction of a suitable Lyapunov function needed in the proof of Lemma~\ref{lm:estb} given below.
\begin{proposition}
\label{prop:est}Let $\mathcal{K}=\diag\left(K_{1},\cdots,K_{N}\right)$ by Lemma~\ref{lm:gain_conv} and let $\mathcal{H}=\diag\left(H_{1},\cdots,H_{N}\right)$. Then, there exist $\Vap_{\mathcal{K}}>0$, a (deterministic) time $t_{\mathcal{K}}$ and a constant $c_{\mathcal{K}}$, such that,
\begin{equation}
\label{prop:est1}
\mathbf{z}^{T}\left(\beta_{t}\overline{L}\otimes I_{M}+\alpha_{t}\widetilde{\mathcal{K}}\mathcal{H}\right)\mathbf{z}\geq c_{\mathcal{K}}\alpha_{t}\left\|\mathbf{z}\right\|^{2},
\end{equation}
for all $t\geq t_{\mathcal{K}}$, $\mathbf{z}\in\mathbb{R}^{NM}$, and $\widetilde{\mathcal{K}}$ satisfying $\|\widetilde{\mathcal{K}}\mathcal{H}-\mathcal{K}\mathcal{H}\|\leq\Vap_{\mathcal{K}}$.
\end{proposition}
\begin{proposition}
\label{prop:cons}
Let $\mathcal{K}$ and $\mathcal{H}$ be defined as in Proposition~\ref{prop:est}. Then, for every $0<\Vap<1$ there exist a deterministic time $t_{\Vap}$ and a constant $c_{\Vap}$, such that,
\begin{equation}
\label{prop:cons1}
\mathbf{z}^{T}\left(\beta_{t}\overline{L}\otimes I_{M}+\alpha_{t}\widetilde{\mathcal{K}}\mathcal{H}\right)\mathbf{z}\geq c_{\Vap}\beta_{t}\left\|\mathbf{z}_{\PC}\right\|^{2}
\end{equation}
for all $t\geq t_{\Vap}$, $\mathbf{z}\in\mathbb{R}^{NM}$ and $\widetilde{\mathcal{K}}$ satisfying
\begin{equation}
\label{prop:cons2}
\left\|\widetilde{\mathcal{K}}\mathcal{H}-\mathcal{K}\mathcal{H}\right\|\leq\Vap.
\end{equation}
Also, in the above $\mathbf{z}_{\PC}$ denotes the projection of $\mathbf{z}$ in the orthogonal complement of the consensus subspace $\mathcal{C}$ as defined in~\eqref{prop:est3} in Appendix~\ref{app:est}.
\end{proposition}

\begin{proof}[Proof of Lemma~\ref{lm:estb}]
The estimator recursions in~\eqref{est_up} may be written as
\begin{equation}
\label{lm:estb10}
\mathbf{x}_{t+1}=\left(I_{NM}-\beta_{t}\overline{L}\otimes I_{M}-\alpha_{t}\mathcal{K}_{t}\mathcal{H}\right)\mathbf{x}_{t}-\beta_{t}\left(\widetilde{L}_{t}\otimes I_{M}\right)\mathbf{x}_{t}+\alpha_{t}\mathcal{K}_{t}\mathbf{y}_{t},
\end{equation}
with $\mathbf{x}_{t}$ and $\mathbf{y}_{t}$ denoting $[\mathbf{x}_{1}^{T}(t),\cdots,\mathbf{x}_{N}^{T}(t)]^{T}$ and $[\mathbf{y}_{1}^{T}(t),\cdots,\mathbf{y}_{N}^{T}(t)]^{T}$ respectively. The sequence $\{\WL_{t}\}$ denotes the sequence of zero mean i.i.d.~ matrices given by $\WL_{t}=L_{t}-\OL_{t}$ for all $t$. The process $\{\mathbf{z}_{t}\}$ defined as $\mathbf{z}_{t}=\mathbf{x}_{t}-\mathbf{1}_{N}\otimes\mathbf{\theta}^{\ast}$ may then be shown to satisfy the recursion
\begin{equation}
\label{lm:estb11}
\mathbf{z}_{t+1}=\left(I_{NM}-\beta_{t}\overline{L}\otimes I_{M}-\alpha_{t}\mathcal{K}_{t}\mathcal{H}\right)\mathbf{z}_{t}-\beta_{t}\left(\widetilde{L}_{t}\otimes I_{M}\right)\mathbf{z}_{t}+\alpha_{t}\mathcal{K}_{t}\mathbf{\zeta}_{t},
\end{equation}
with $\mathbf{\zeta}_{t}=[\mathbf{\zeta}_{1}^{T}(t),\cdots,\mathbf{\zeta}_{N}^{T}(t)]^{T}$. Now fix $0<\Vap<\Vap_{\mathcal{K}}\wedge 1$, where $\Vap_{\mathcal{K}}$ is defined in the hypothesis of Proposition~\ref{prop:est}. Since, $\mathcal{K}_{t}\rightarrow\mathcal{K}$ a.s., by Egorov's theorem (\cite{Halmos}) for every $\delta>0$, there exists $t_{\delta}$ such that
\begin{equation}
\label{lm:estb12}
\Past\left(\sup_{t\geq t_{\delta}}\left\|\mathcal{K}_{t}\mathcal{H}-\mathcal{K}\mathcal{H}\right\|\leq\Vap\right)>1-\delta~~\mbox{and}~~\Past\left(\sup_{t\geq t_{\delta}}\left\|\mathcal{K}_{t}-\mathcal{K}\right\|\leq\Vap\right)>1-\delta.
\end{equation}
Moreover, such a $t_{\delta}$ may be chosen to satisfy $t_{\delta}>t_{\mathcal{K}}\vee t_{\Vap}$, where $t_{\mathcal{K}}$ and $t_{\Vap}$ are defined in the hypotheses of Proposition~\ref{prop:est} and Proposition~\ref{prop:cons}, respectively.

Let $\mathcal{K}_{\Vap}$ be a (deterministic) matrix, such that,
\begin{equation}
\label{lm:estb14}
\left\|\mathcal{K}_{\Vap}\mathcal{H}-\mathcal{K}\mathcal{H}\right\|<\Vap~~~\mbox{and}~~~\left\|\mathcal{K}_{\Vap}-\mathcal{K}\right\|<\Vap.
\end{equation}
Then, for every $\delta>0$, we may define the $\{\mathcal{F}_{t}\}$ adapted process $\{\mathcal{K}_{t}^{\delta}\}$, such that,
\begin{equation}
\label{lm:estb15}
\mathcal{K}_{t}^{\delta}=\left\{\begin{array}{ll}
                                \mathcal{K}_{t} & \mbox{if $t<t_{\delta}$}\\
                                \mathcal{K}_{t} & \mbox{if $t\geq t_{\delta}$ and $\|\mathcal{K}_{t}\mathcal{H}-\mathcal{K}\mathcal{H}\|\vee\|\mathcal{K}_{t}-\mathcal{K}\|\leq\Vap$}\\
                                \mathcal{K}_{\Vap} & \mbox{otherwise}.
                                \end{array}
                                \right.
\end{equation}
Also, for each $\delta>0$, we define the $\{\mathcal{F}_{t}\}$ adapted process $\{\mathbf{z}_{t}^{\delta}\}$ by the recursion
\begin{equation}
\label{lm:estb16}
\mathbf{z}_{t+1}^{\delta}=\left(I_{NM}-\beta_{t}\overline{L}\otimes I_{M}-\alpha_{t}\mathcal{K}_{t}^{\delta}\mathcal{H}\right)\mathbf{z}_{t}^{\delta}-\beta_{t}\left(\widetilde{L}_{t}\otimes I_{M}\right)\mathbf{z}_{t}^{\delta}+\alpha_{t}\mathcal{K}_{t}^{\delta}\mathbf{\zeta}_{t},
\end{equation}
with $\mathbf{z}_{0}^{\delta}=\mathbf{z}_{0}$. To show that the process $\{\mathbf{z}_{t}\}$ (and, hence $\{\mathbf{x}_{t}\}$) is bounded a.s., we note that it suffices to show that the process $\{\mathbf{z}_{t}^{\delta}\}$ is bounded a.s. for each $\delta>0$. This is due to the fact that, by the definition of $t_{\delta}$, for each $\delta>0$ we have
\begin{equation}
\label{lm:estb17}
\Past\left(\sup_{t\geq 0}\left\|\mathcal{K}_{t}^{\delta}-\mathcal{K}_{t}\right\|=0\right)>1-\delta,
\end{equation}
and, hence
\begin{equation}
\label{lm:estb18}
\Past\left(\sup_{t\geq 0}\left\|\mathbf{z}_{t}^{\delta}-\mathbf{z}_{t}\right\|=0\right)>1-\delta.
\end{equation}
Thus the boundedness of the processes $\{\mathbf{z}_{t}^{\delta}\}$ for each $\delta>0$ would imply
\begin{equation}
\label{lm:estb19}
\Past\left(\sup_{t\geq 0}\left\|\mathbf{x}_{t}\right\|<\infty\right)>1-\delta
\end{equation}
for every $\delta>0$. The assertion of Lemma~\ref{lm:estb} would then follow by taking $\delta$ to zero.

Hence, in the following, we focus only on the processes $\{\mathbf{z}^{\delta}_{t}\}$ and show that the latter are bounded a.s. for every $\delta>0$. To this end, fix $\delta>0$ and consider the $\mathcal{F}_{t}$ process $V_{t}^{\delta}=\|\mathbf{z}_{t}^{\delta}\|^{2}$. It can be shown (Assumption~\textbf{(A.3)}) that
\begin{eqnarray}
\label{lm:estb20}
\East\left[V_{t+1}^{\delta}~|~\mathcal{F}_{t}\right] & = & V_{t}^{\delta}+\beta_{t}^{2}\left(\mathbf{z}_{t}^{\delta}\right)^{T}\East\left[(\WL_{t}\otimes I_{M})^{2}\right]\mathbf{z}_{t}^{\delta}+\alpha_{t}^{2}\East\left[\left\|\mathcal{K}_{t}^{\delta}\mathbf{\zeta}_{t}\right\|^{2}\right]\nonumber \\ & & -2\left(\mathbf{z}_{t}^{\delta}\right)^{T}\left(\beta_{t}\OL\otimes I_{M}+\alpha_{t}\mathcal{K}_{t}^{\delta}\mathcal{H}\right)\mathbf{z}_{t}^{\delta}+\beta_{t}^{2}\left(\mathbf{z}_{t}^{\delta}\right)^{T}\left(\OL\otimes I_{M}\right)^{2}\mathbf{z}_{t}^{\delta}\nonumber \\ & & + \alpha_{t}^{2}\left(\mathbf{z}_{t}^{\delta}\right)^{T}\left(\mathcal{K}_{t}^{\delta}\mathcal{H}\right)^{T}\mathcal{K}_{t}^{\delta}\mathcal{H}\mathbf{z}_{t}^{\delta}+2\alpha_{t}\beta_{t}\left(\mathbf{z}_{t}^{\delta}\right)^{T}\left(\OL\otimes I_{M}\right)\left(\mathcal{K}_{t}^{\delta}\mathcal{H}\right)\mathbf{z}_{t}^{\delta}.
\end{eqnarray}
Since the Laplacians are bounded matrices by definition and the matrix $\mathcal{K}_{t}^{\delta}$ is bounded for $t\geq t_{\delta}$ by construction, there exists a constant $c_{3}>0$, sufficiently large, such that the inequalities
\begin{equation}
\label{lm:estb21}
\left(\mathbf{z}_{t}^{\delta}\right)^{T}\East\left[(\WL_{t}\otimes I_{M})^{2}\right]\mathbf{z}_{t}^{\delta}=\left(\mathbf{z}_{t,\PC}^{\delta}\right)^{T}\East\left[(\WL_{t}\otimes I_{M})^{2}\right]\mathbf{z}_{t,\PC}^{\delta}\leq c_{3}\left\|\mathbf{z}_{t,\PC}^{\delta}\right\|^{2},
\end{equation}
\begin{equation}
\label{lm:estb22}
\left(\mathbf{z}_{t}^{\delta}\right)^{T}\left(\OL\otimes I_{M}\right)^{2}\mathbf{z}_{t}^{\delta}=\left(\mathbf{z}_{t,\PC}^{\delta}\right)^{T}\left(\OL\otimes I_{M}\right)^{2}\mathbf{z}_{t,\PC}^{\delta}\leq c_{3}\left\|\mathbf{z}_{t,\PC}^{\delta}\right\|^{2},
\end{equation}
\begin{equation}
\label{lm:estb23}
\left(\mathbf{z}_{t}^{\delta}\right)^{T}\left(\OL\otimes I_{M}\right)\left(\mathcal{K}_{t}^{\delta}\mathcal{H}\right)\mathbf{z}_{t}^{\delta}
\leq c_{3}\left\|\mathbf{z}_{t}^{\delta}\right\|^{2},
\end{equation}
\begin{equation}
\label{lm:estb24}
\left(\mathbf{z}_{t}^{\delta}\right)^{T}\left(\mathcal{K}_{t}^{\delta}\mathcal{H}\right)^{T}\mathcal{K}_{t}^{\delta}\mathcal{H}\mathbf{z}_{t}^{\delta}
\leq c_{3}\left\|\mathbf{z}_{t}^{\delta}\right\|^{2},~~
\East\left[\left\|\mathcal{K}_{t}^{\delta}\mathbf{\zeta}_{t}\right\|^{2}\right]\leq c_{3}
\end{equation}
hold for all $t\geq t_{\delta}$ with $\mathbf{z}_{t,\PC}^{\delta}$ denoting the projection of $\mathbf{z}_{t}^{\delta}$ on the subspace $\PC$. Also, by Proposition~\ref{prop:est} and Proposition~\ref{prop:cons}, for $t\geq t_{\delta}$,
\begin{equation}
\label{lm:estb2}
\left(\mathbf{z}_{t}^{\delta}\right)^{T}\left(\beta_{t}\OL\otimes I_{M}+\alpha_{t}\mathcal{K}_{t}^{\delta}\mathcal{H}\right)\mathbf{z}_{t}^{\delta}\geq c_{\mathcal{K}}\alpha_{t}\left\|\mathbf{z}_{t}^{\delta}\right\|^{2}+c_{\Vap}\beta_{t}\left\|\mathbf{z}_{t,\PC}^{\delta}\right\|^{2},
\end{equation}
where the positive constants $c_{\mathcal{K}}$ and $c_{\Vap}$ are defined in the hypotheses of Proposition~\ref{prop:est} and Proposition~\ref{prop:cons}, respectively. The inequalities~\eqref{lm:estb21}-\eqref{lm:estb24} then lead to
\begin{eqnarray}
\label{lm:estb27}
\East\left[V_{t+1}^{\delta}~|~\mathcal{F}_{t}\right] & \leq & V_{t}-\left(c_{\mathcal{K}}\beta_{t}-2c_{3}\beta_{t}^{2}\right)\left\|\mathbf{z}_{t,\PC}^{\delta}\right\|^{2}\nonumber 
\end{eqnarray}
for all $t\geq t_{\delta}$. Observing the decay rates of the various terms in~\eqref{weight}, we conclude that there exists $\bar{t}_{\delta}\geq t_{\delta}$, such that,
\begin{equation}
\label{lm:estb28}
c_{\mathcal{K}}\beta_{t}-2c_{3}\beta_{t}^{2}>0~~\mbox{and}~~c_{\Vap}\alpha_{t}-2\alpha_{t}\beta_{t}c_{3}-\alpha_{t}^{2}c_{3}>0,
\end{equation}
for $t\geq\bar{t}_{\delta}$ and, hence,
\begin{equation}
\label{lm:estb29}
\East\left[V_{t+1}^{\delta}~|~\mathcal{F}_{t}\right]\leq V_{t}^{\delta}+c_{3}\alpha_{t}^{2}
\end{equation}
for all $t\geq\bar{t}_{\delta}$. Let us introduce the $\{\mathcal{F}_{t}\}$ adapted process $\{\overline{V}^{\delta}_{t}\}$, such that,
\begin{equation}
\label{lm:estb30}
\overline{V}^{\delta}_{t}=V_{t}^{\delta}-c_{3}\sum_{s=t}^{\infty}\alpha_{s}^{2}
\end{equation}
for $t\geq 0$. The process $\{\overline{V}^{\delta}_{t}\}$ is well-defined as the sequence $\{\alpha_{t}\}$ is square summable. From~\eqref{lm:estb29} it follows immediately that
\begin{equation}
\label{lm:estb31}
\East\left[\overline{V}^{\delta}_{t+1}~|~\mathcal{F}_{t}\right]\leq V_{t}^{\delta}-c_{3}\alpha_{t}^{2}-c_{3}\sum_{s=t+1}^{\infty}\alpha_{s}^{2}=\overline{V}^{\delta}_{t}
\end{equation}
for $t\geq\overline{t}_{\delta}$. Hence, the process $\{\overline{V}^{\delta}_{t}\}_{t\geq\overline{t}_{\delta}}$ is a supermartingale. Moreover, it is bounded from below, since $V_{t}\geq 0$ by construction, and, in fact,
\begin{equation}
\label{lm:estb32}
\overline{V}^{\delta}_{t}\geq -c_{3}\sum_{s=0}^{\infty}\alpha_{s}^{2}
\end{equation}
for all $t\geq 0$. Thus $\{\overline{V}^{\delta}_{t}\}_{t\geq\overline{t}_{\delta}}$ is a supermartingale that is bounded from below and, hence converges a.s. to a finite random variable $\overline{V}^{\delta}$, i.e., $\overline{V}_{t}^{\delta}\rightarrow\overline{V}$ a.s. as $t\rightarrow\infty$. In particular, the process $\{\overline{V}_{t}^{\delta}\}$ is pathwise bounded. By~\eqref{lm:estb30} the process $\{V_{t}^{\delta}\}$ is also pathwise bounded. Thus, for each $\delta>0$, the process $\{\mathbf{z}_{t}^{\delta}\}$ is bounded a.s. and the assertion follows.
\end{proof}

The next result (see Appendix~\ref{app:est} for a proof) quantifies the rate at which the different agent estimates reach agreement and is stated as follows:
\begin{lemma}
\label{lm:estc}
For every $\tau_{0}$ such that
$0\leq\tau_{0}<\tau_{1}-\tau_{2}-1/(2+\Vap_{1})$, we have
\begin{equation}
\label{lm:estc2}
\Past\left(\lim_{t\rightarrow\infty}(t+1)^{\tau_{0}}\left(\mathbf{x}_{n}(t)-\Xavg(t)\right)=0\right)=1
\end{equation}
with $\Xavg(t)=(1/N)\sum_{n=1}^{N}\mathbf{x}_{n}(t)$ denoting the instantaneous network averaged estimate.
\end{lemma}

The rest of the section focuses on the convergence properties of the network averaged estimate $\{\Xavg(t)\}$ and completes the final steps required to establish the convergence properties of the agent estimates $\{\mathbf{x}_{n}(t)\}$. The first result in this direction concerns the consistency of the average estimate sequence.
\begin{lemma}
\label{lm:avgc}Under the additional assumption that $\tau_{1}=1$ (see~\textbf{(A.5)}) we have
\begin{equation}
\label{lm:avgc1}
\Past\left(\lim_{t\rightarrow\infty}\left(\Xavg(t)-\mathbf{\theta}^{\ast}\right)=0\right)=1
\end{equation}
with $\Xavg(t)=(1/N)\sum_{n=1}^{N}\mathbf{x}_{n}(t)$ the instantaneous network averaged estimate.
\end{lemma}
\begin{proof}
Let us denote by $\mathbf{z}_{t}$ the residual $\Xavg(t)-\mathbf{\theta}^{\ast}$. The $\mathcal{F}_{t}$-adapted process $\{\mathbf{z}_{t}\}$ may be shown to satisfy the recursion
\begin{equation}
\label{lm:avgc2}
\mathbf{z}_{t+1}=\left(I_{M}-\alpha_{t}\Gamma_{t}\right)\mathbf{z}_{t}+\alpha_{t}U_{t}+\alpha_{t}J_{t}
\end{equation}
with $\{\Gamma_{t}\}$, $\{U_{t}\}$ being $\mathcal{F}_{t}$-adapted, $\{J_{t}\}$ being $\mathcal{F}_{t+1}$-adapted and given by
\begin{equation}
\label{lm:avgc3}
\Gamma_{t}=\frac{1}{N}\sum_{n=1}^{N}K_{n}(t)H_{n},~~U_{t}=\frac{1}{N}\sum_{n=1}^{N}K_{n}(t)\left(\mathbf{x}_{n}(t)-\Xavg(t)\right)~\mbox{and}~J_{t}=\frac{1}{N}K_{n}(t)\mathbf{\zeta}_{n}(t)
\end{equation}
respectively. Now fix $0<\tau_{0}<\tau_{1}-\tau_{2}-1/(2+\Vap_{1})$ and, by the convergence of the gain processes and Lemma~\ref{lm:estc}, $\Gamma_{t}\rightarrow I_{M}$ and $(t+1)^{\tau_{0}}U_{t}\rightarrow 0$ a.s. as $t\rightarrow\infty$. By Egorov's theorem the a.s. convergence may be assumed to be uniform on sets of arbitrarily large probability measure and, hence, for every $\delta>0$, there exist uniformly bounded processes $\{\Gamma_{t}^{\delta}\}$, $\{U_{t}^{\delta}\}$, and $\{\mathcal{K}_{t}^{\delta}\}$ satisfying
\begin{equation}
\label{lm:avgc5}
\Past\left(\sup_{s\geq t_{\Vap}^{\delta}}\left\|\Gamma_{s}^{\delta}-I_{M}\right\|\vee\left\|\mathcal{K}^{\delta}_{t}-\mathcal{K}\right\|>\Vap\right)=0~\mbox{and}~
\Past\left(\sup_{s\geq t_{\Vap}^{\delta}}(s+1)^{\tau_{0}}\left\|U_{s}^{\delta}\right\|>\Vap\right)=0
\end{equation}
for each $\Vap>0$ and some $t_{\Vap}^{\delta}$ (sufficiently large), such that
\begin{equation}
\label{lm:avgc7}
\Past\left(\sup_{t\geq 0}\left\|\Gamma_{t}^{\delta}-\Gamma_{t}\right\|\vee\left\|\mathcal{K}^{\delta}_{t}-\mathcal{K}_{t}\right\|\vee\left\|U_{t}^{\delta}-U_{t}\right\|=0\right)>1-\delta.
\end{equation}
Also, for each $\delta>0$, define the $\mathcal{F}_{t}$-adapted process $\{\mathbf{z}^{\delta}_{t}\}$ by
\begin{equation}
\label{lm:avgc8}
\mathbf{z}_{t+1}^{\delta}=\left(I_{M}-\alpha_{t}\Gamma_{t}^{\delta}\right)\mathbf{z}_{t}^{\delta}+\alpha_{t}U_{t}^{\delta}+\alpha_{t}J_{t}^{\delta}
\end{equation}
with $\mathbf{z}_{0}^{\delta}=\mathbf{z}_{0}$ and $J_{t}^{\delta}=\frac{1}{N}\sum_{n=1}^{N}K_{n}^{\delta}(t)\mathbf{\zeta}_{n}(t)$ and
\begin{equation}
\label{lm:avgc9}
\Past\left(\sup_{t\geq 0}\left\|\mathbf{z}_{t}^{\delta}-\mathbf{z}_{t}\right\|=0\right)>1-\delta.
\end{equation}
By the above development, to show that $\mathbf{z}_{t}\rightarrow 0$ as $t\rightarrow\infty$, it suffices to show that $\mathbf{z}_{t}^{\delta}\rightarrow 0$ as $t\rightarrow\infty$ for each $\delta>0$. Hence, in the following, we focus on the process $\{\mathbf{z}_{t}^{\delta}\}$ only for a fixed but arbitrary $\delta>0$.

Now let $\{V_{t}^{\delta}\}$ denote the $\{\mathcal{F}_{t}\}$ adapted process such that $V_{t}^{\delta}=\left\|\mathbf{z}_{t}^{\delta}\right\|^{2}$ for all $t$. Using the fact that $\East\left[\mathbf{\zeta}_{t}~|~\mathcal{F}_{t}\right]=\mathbf{0}$ for all $t$, it follows that
\begin{eqnarray}
\label{lm:avgc10}
\East\left[V_{t+1}^{\delta}~|~\mathcal{F}_{t}\right] &\leq & \left\|I_{M}-\alpha_{t}\Gamma_{t}^{\delta}\right\|^{2}V_{t}^{\delta}+2\alpha_{t}(U_{t}^{\delta})^{T}\left(I_{M}-\alpha_{t}\Gamma_{t}^{\delta}\right)\mathbf{z}_{t}^{\delta}\nonumber \\ & & +\alpha_{t}^{2}\left\|U_{t}\right\|^{2}+\alpha_{t}^{2}\East\left[\left\|J_{t}\right\|^{2}~|~\mathcal{F}_{t}\right].
\end{eqnarray}
For $t$ large enough
\begin{equation}
\label{lm:avgc11}
\left|2\alpha_{t}U_{t}^{T}\left(I_{M}-\alpha_{t}\Gamma_{t}^{\delta}\right)\mathbf{z}_{t}^{\delta}\right|\leq 2\alpha_{t}\left\|U_{t}^{\delta}\right\|\left\|\mathbf{z}_{t}^{\delta}\right\|\leq 2\alpha_{t}\left\|U_{t}^{\delta}\right\|\left\|\mathbf{z}_{t}^{\delta}\right\|^{2}+2\alpha_{t}\left\|U_{t}^{\delta}\right\|.
\end{equation}
Then making $t_{\Vap}^{\delta}$ larger (if necessary), such that $\|U_{t}^{\delta}\|\leq\Vap (t+1)^{-\tau_{0}}$, $\East[\|J_{t}\|^{2}|\mathcal{F}_{t}]$ is uniformly bounded, and~\eqref{lm:avgc11} holds for all $t\geq t_{\Vap}^{\delta}$, it follows from~\eqref{lm:avgc10}-\eqref{lm:avgc11} that there exist positive constants $c_{1}$ and $c_{2}$ so that
\begin{align}
\label{lm:avgc12}
\East\left[V_{t+1}^{\delta}~|~\mathcal{F}_{t}\right]\leq&
\left(1-c_{1}\alpha_{t}+c_{2}\alpha_{t}(t+1)^{-\tau_{0}}\right)V_{t}^{\delta}\\
+&c_{2}\left(\alpha_{t}(t+1)^{-\tau_{0}}+\alpha_{t}^{2}(t+1)^{-2\tau_{0}}+\alpha_{t}^{2}\right)
\end{align}
for all $t\geq t_{\Vap}^{\delta}$. Since $0<\tau_{0}<\tau_{1}$, the first term inside the second set of parenthesis on the right hand side dominates; by making $c_{4}>c_{2}$ and $c_{3}<c_{1}$ appropriately, we have
\begin{equation}
\label{lm:avgc13}
\East\left[V_{t+1}^{\delta}~|~\mathcal{F}_{t}\right]\leq\left(1-c_{3}\alpha_{t}\right)V_{t}^{\delta}+c_{4}\alpha_{t}(t+1)^{-\tau_{0}}\leq V_{t}^{\delta}+c_{4}\alpha_{t}(t+1)^{-\tau_{0}}
\end{equation}
for all $t\geq t_{\Vap}^{\delta}$. Now consider the $\{\mathcal{F}_{t}\}$ adapted process $\{\overline{V}_{t}^{\delta}\}$, such that,
\begin{equation}
\label{lm:avgc14}
\overline{V}_{t}^{\delta}=V_{t}^{\delta}-c_{4}\sum_{s=t}^{\infty}\alpha_{s}(s+1)^{-\tau_{0}}
\end{equation}
for $t\geq 0$. Since $\tau_{1}=1$ and $\tau_{0}>0$, the sequence $\{\alpha_{t}(t+1)^{-\tau_{0}}\}$ is summable and the process $\{\overline{V}_{t}^{\delta}\}$ is bounded from below. It is readily seen that $\{\overline{V}_{t}^{\delta}\}_{t\geq t^{\delta}_{\Vap}}$ is a supermartingale and, hence converges a.s. to a finite random variable. By~\eqref{lm:avgc14}, the process $\{V_{t}^{\delta}\}$ also converges a.s. to a finite random variable $V^{\delta}$ (necessarily non-negative). Finally, from~\eqref{lm:avgc13},
\begin{equation}
\label{lm:avgc15}
\East\left[V_{t+1}^{\delta}\right]\leq\left(1-c_{3}\alpha_{t}\right)\East\left[V_{t}^{\delta}\right]+c_{4}\alpha_{t}(t+1)^{-\tau_{0}}
\end{equation}
for $t\geq t_{\Vap}^{\delta}$. Since $\tau_{0}>0$ the sequence $\{\alpha_{t}(t+1)^{-\tau_{0}}\}$ decays faster than $\{\alpha_{t}\}$ and, hence by Lemma~\ref{lm:JSTSP-det} we have $\East[V_{t}^{\delta}]\rightarrow 0$ as $t\rightarrow\infty$. The sequence $\{V_{t}^{\delta}\}$ is non-negative, so by Fatou's lemma we further conclude that
\begin{equation}
\label{lm:avgc16}
0\leq\East\left[V^{\delta}\right]\leq\liminf_{t\rightarrow\infty}\East\left[V_{t}^{\delta}\right]=0.
\end{equation}
The above implies $V^{\delta}=0$ a.s. by the non-negativity of $V^{\delta}$. Hence $\|\mathbf{z}_{t}^{\delta}\|\rightarrow 0$ as $t\rightarrow\infty$ and the desired assertion follows.
\end{proof}

By inductive reasoning, we now obtain a stronger version of Lemma~\ref{lm:avgc} that quantifies the convergence rate in the above (see Appendix~\ref{app:est} for a proof).
\begin{lemma}
\label{lm:avgr} Let assumptions~\textbf{(A.1)}-\textbf{(A.5)} hold with $\tau_{1}=1$ and $a\geq 1$. Then, for each $n$ and $\tau\in [0,1/2)$,
\begin{equation}
\label{lm:avgr1}
\Past\left(\lim_{t\rightarrow\infty}(t+1)^{\tau}\left\|\mathbf{x}_{n}(t)-\mathbf{\theta}^{\ast}\right\|=0\right)=1.
\end{equation}
\end{lemma}

\section{Proofs of Main Results}
\label{sec:proof_main_res} The proof of Theorem~\ref{th:estc} is a direct consequence of the triangle inequality and Lemma~\ref{lm:estc} since all agent estimates converge to the network-averaged estimate at the required rate.

\textbf{Proof of Theorem~\ref{th:estcons}}

\begin{proof}
Since $\Vap_{1}>0$, $\tau_{1}=1$ and $\tau_{1}>\tau_{2}+1/(2+\Vap_{1})+1/2$, from Lemma~\ref{lm:estc}  there exists $\Vap>0$ (sufficiently small) such that
\begin{equation}
\label{th:estcons10}
\Past\left(\lim_{t\rightarrow\infty}(t+1)^{1/2+\Vap}\left\|\mathbf{x}_{n}(t)-\Xavg(t)\right\|=0\right)=1
\end{equation}
for all $n$. Moreover, by Lemma~\ref{lm:avgr}, for each $\tau\in [0,1/2)$, we have $(t+1)^{\tau}\|\Xavg(t)-\mathbf{\theta}^{\ast}\|\rightarrow 0$ a.s. as $t\rightarrow\infty$, for all $n$. Since $\tau<1/2+\Vap$, an immediate application of the triangle inequality yields the required estimate convergence rate.
\end{proof}

\textbf{Proof of Theorem~\ref{th:estn}}

We will use the following result from~\cite{Fabian-2} concerning the asymptotic normality of non-Markov stochastic recursions. The statement here is somewhat less general than in~\cite{Fabian-2} but serves our application and eases the additional notational complexity.
\begin{lemma}[Theorem 2.2.~in~\cite{Fabian-2}]
\label{lm:Fab-2} Let $\{\mathbf{z}_{t}\}$ be an $\mathbb{R}^{k}$ valued $\{\mathcal{F}_{t}\}$ adapted process that satisfies
\begin{equation}
\label{lm:Fab-21}
\mathbf{z}_{t+1}=\left(I_{k}-\frac{1}{t+1}\Gamma_{t}\right)\mathbf{z}_{t}+(t+1)^{-1}\Phi_{t}V_{t}+(t+1)^{-3/2}T_{t},
\end{equation}
where $\{V_{t}\}$ and $\{T_{t}\}$ are $\mathbb{R}^{k}$ valued stochastic processes, such that, for each $t$, $V_{t-1}$ and $T_{t}$ are $\mathcal{F}_{t}$-adapted, and where the processes $\{\Gamma_{t}\}$ and $\{\Phi_{t}\}$ are $\mathbb{R}^{k\times k}$ valued and $\{\mathcal{F}_{t}\}$ adapted. Assume
\begin{equation}
\label{lm:Fab-22}
\Gamma_{t}\rightarrow I_{k},~~\Phi_{t}\rightarrow\Phi~~\mbox{and}~~T_{t}\rightarrow\mathbf{0}~~\mbox{as $t\rightarrow\infty$.}
\end{equation}
Let the sequence $\{V_{t}\}$ satisfy $\mathbb{E}[V_{t}|\mathcal{F}_{t}]=\mathbf{0}$ for each $t$ and there exist a constant $C>0$ and a matrix $\Sigma$ such that $C>\left\|\mathbb{E}[V_{t}V_{t}^{T}|\mathcal{F}_{t}]-\Sigma\right\|\rightarrow 0$ as $t\rightarrow\infty$, and, with
\begin{equation}
\label{lm:Fab-24}
\sigma_{t,r}^{2}=\int_{\|V_{t}\|^{2}\geq r(t+1)}\|V_{t}\|^{2}d\mathbb{P},
\end{equation}
let $\lim_{t\rightarrow\infty}\frac{1}{t+1}\sum_{s=0}^{t}\sigma^{2}_{s,r}=0$ for every $r>0$. Then, the asymptotic distribution of $(t+1)^{-1/2}\mathbf{z}_{t}$ is normal with mean $\mathbf{0}$ and covariance matrix $\Phi\Sigma\Phi^{T}$.
\end{lemma}

\begin{proof}[Proof of Theorem~\ref{th:estn}]
The residual process $\{\mathbf{z}_{t}\}$ and its $\delta$-approximations $\{\mathbf{z}_{t}^{\delta}\}$ are given in~\eqref{lm:avgc2}-\eqref{lm:avgc8}. With $\tau_{1}=a=1$,
\begin{equation}
\label{th:estn10}
\mathbf{z}_{t+1}=\left(I_{M}-\frac{1}{t+1}\Gamma_{t}\right)\mathbf{z}_{t}+(t+1)^{-1}U_{t}+(t+1)^{-1}J_{t},
\end{equation}
where $U_{t}$ and $J_{t}$ are defined in~\eqref{lm:avgc2}-\eqref{lm:avgc8}. Since $J_{t}=(1/N)\sum_{n=1}^{N}K_{n}(t)\mathbf{\zeta}_{n}(t)$ and the $\{K_{n}(t)\}$'s may not converge uniformly (both in time and space), Lemma~\ref{lm:Fab-2} is not applicable directly. Hence, we first consider the process $\{\mathbf{z}_{t}^{\delta}\}$ for some $\delta>0$. In order to apply Lemma~\ref{lm:Fab-2} to the process $\{\mathbf{z}_{t}^{\delta}\}$, define
\begin{equation}
\label{th:estn11}
T_{t}=(t+1)^{1/2}U_{t}^{\delta}
\end{equation}
for each $t$. Note that by~\eqref{lm:avgr4} $\|U_{t}^{\delta}\|=o\left((t+1)^{-1/2}\right)$ and, hence $T_{t}\rightarrow\mathbf{0}$ as $t\rightarrow\infty$. Also define
\begin{equation}
\Phi_{t}=I_{M}~~\mbox{and}~~V_{t}=J_{t}^{\delta}
\end{equation}
for each $t$. Clearly, $\East[V_{t}|\mathcal{F}_{t}]=\mathbf{0}$ for all $t$. By the convergence of $\mathcal{K}_{t}^{\delta}$ to $\mathcal{K}$,
\begin{equation}
\label{th:estn12}
\lim_{t\rightarrow\infty}\East\left[V_{t}V_{t}^{T}~|~\mathcal{F}_{t}\right] = \lim_{t\rightarrow\infty}\frac{1}{N^{2}}\sum_{n=1}^{N}K_{n}^{\delta}(t)R_{n}\left(K_{n}^{\delta}\right)^{T}=\Sigmac^{-1},
\end{equation}
where the last step follows from Lemma~\ref{lm:gain_conv}. Moreover the uniform boundedness of the process $\{\mathcal{K}_{t}^{\delta}\}$ implies the existence of a constant $C>0$ such that
\begin{equation}
\label{th:estn13}
\left\|\East\left[V_{t}V_{t}^{T}~|~\mathcal{F}_{t}\right]-\Sigmac^{-1}\right\|<C
\end{equation}
for all $t\geq 0$. The $\{V_{t}\}$ thus constructed also satisfies the uniform integrability assumption~\eqref{lm:Fab-24} due to the i.i.d. nature of the noise processes and the uniform boundedness of $\{\mathcal{K}_{t}^{\delta}\}$. Thus, the process $\{\mathbf{z}_{t}^{\delta}\}$ falls under the purview of Lemma~\ref{lm:Fab-2} with $\Phi=I_{M}$ and $\Sigma=\Sigmac^{-1}$. We thus conclude that
\begin{equation}
\label{th:estn14}
(t+1)^{-1/2}\mathbf{z}_{t}^{\delta}\Longrightarrow\mathcal{N}\left(\mathbf{0},\Sigmac^{-1}\right)
\end{equation}
for each $\delta>0$. To extend this asymptotic normality to the process $\{\mathbf{z}_{t}\}$, consider any bounded continuous function $f:\mathbb{R}^{M}\longmapsto\mathbb{R}$. By weak convergence (Portmanteau's theorem,~\cite{Billingsley}) we have
\begin{equation}
\label{th:estn15}
\lim_{t\rightarrow\infty}\East\left[f\left((t+1)^{-1/2}\mathbf{z}_{t}^{\delta}\right)\right]=\East\left[f\left(\mathbf{z}^{\ast}\right)\right]
\end{equation}
for each $\delta$, where $\mathbf{z}^{\ast}$ denotes a $\mathcal{N}\left(\mathbf{0},\Sigmac^{-1}\right)$ distributed random vector under the measure $\mathbb{P}_{\ast}$. Denoting by $\|f\|_{\infty}$ the sup-norm of $f(\cdot)$ (necessarily finite) we obtain from~\eqref{lm:avgc9}
\begin{equation}
\label{th:estn16}
\left\|\East\left[f\left((t+1)^{-1/2}\mathbf{z}_{t}^{\delta}\right)\right]-\East\left[f\left((t+1)^{-1/2}\mathbf{z}_{t}\right)\right]\right\|\leq 2\delta\|f\|_{\infty}.
\end{equation}
By~\eqref{th:estn15} we then have
\begin{equation}
\label{th:estn17}
\limsup_{t\rightarrow\infty}\left\|\East\left[f\left((t+1)^{-1/2}\mathbf{z}_{t}\right)\right]-\East\left[f\left(\mathbf{z}^{\ast}\right)\right]\right\|\leq 2\delta\|f\|_{\infty}.
\end{equation}
Since the above holds for each $\delta>0$, we conclude that $\East\left[f\left((t+1)^{-1/2}\mathbf{z}_{t}\right)\right]\rightarrow\East\left[f\left(\mathbf{z}^{\ast}\right)\right]$ as $t\rightarrow\infty$. This convergence holds for all bounded continuous functions $f(\cdot)$ thus giving the required weak convergence of the sequence $\left\{(t+1)^{-1/2}\mathbf{z}_{t}\right\}$.
\end{proof}

\section{Conclusion}
\label{conclusion} We have developed a distributed estimator that combines a recursive collaborative learning step with the estimate update task. Through this learning process, the agents adaptively improve their quantitative model information and innovation gains with a view toward achieving the performance of the optimal centralized estimator. Intuitively, the distributed approach is a culmination of two potentials, the agreement (or consensus) and the innovation. By properly designing the relative strength of their excitations, we have shown that the agent estimates may be made asymptotically efficient in terms of their asymptotic covariance which coincides with the asymptotic covariance (the inverse of the Fisher information rate for Gaussian systems) of a centralized estimator with perfect statistical information and having access to all agent observations at all times. A typical application scenario involves multi-sensor distributed platforms, for example, the smart grid or vehicular networks. Such networks are generally equipped with rich sensing infrastructures and high sensing diversity, but suffer from lack of information about the global model and about the relative observation efficiencies due to unpredictable changes and constraints in the sensing resources. Extensions of this work to nonlinear sensing platforms are currently being investigated. Another important direction will be the extension of this adaptive collaborative scheme to dynamic parameter situations as opposed to the static case considered in this paper.\\

\appendix

\section{Proofs in Section~\ref{app_res}}
\label{app:app_res}
{\small
\begin{proof}{Proof of Lemma~\ref{lm:seq-gen}}
We start by showing that for each positive integer $k$, the following holds:
\begin{equation}
\label{lm:seq1}
\lim_{t\rightarrow\infty}(t+1)^{k(\delta_{2}-\delta_{1}-\varepsilon_{0})}\mathbb{E}\left[\mathbf{z}_{t}^{k}\right]=0
\end{equation}
for every $0<\varepsilon_{0}\leq\delta_{2}-\delta_{1}$. The proof proceeds by induction on $k$. Let us first consider $k=1$. We then have
\begin{eqnarray}
\label{lm:seq2}
\mathbb{E}\left[\mathbf{z}_{t+1}\right] & \leq & \mathbb{E}\left[\left(1-\mathbb{E}[r_{1}(t)~|~\mathcal{F}_{t}]\right)\mathbf{z}_{t}\right]+r_{2}(t)\nonumber \\ & \leq & (1-\overline{r}_{1}(t))\mathbb{E}[\mathbf{z}_{t}]+r_{2}(t),
\end{eqnarray}
where by $\overline{r}_{1}(t)$ we denote the quantity $a_{1}/(t+1)^{\delta_{1}}$. The deterministic $\mathbb{R}_{+}$ valued sequence $\{\mathbb{E}[\mathbf{z}_{t}]\}$ satisfies the conditions of Lemma~\ref{lm:JSTSP-det} and the claim in~\eqref{lm:seq1} holds for $k=1$.
Now assume the claim in~\eqref{lm:seq1} holds for all $k\leq k_{0}$, with $k_{0}$ a positive integer. We now show that the claim also holds for $k=k_{0}+1$. Indeed, by the polynomial expansion
\begin{equation}
\label{lm:seq3}
\mathbf{z}_{t+1}^{k_{0}+1}=\sum_{i=0}^{k_{0}+1}{k_{0}+1\choose i}\left(\left(1-r_{1}(t)\right)\mathbf{z}_{t}\right)^{k_{0}+1-i}r_{2}^{i}(t)
\end{equation}
and the fact that $0\leq r_{1}(t)\leq 1$, we have
\begin{equation}
\label{lm:seq4}
\mathbf{z}_{t+1}^{k_{0}+1}\leq (1-r_{1}(t))\mathbf{z}_{t}^{k_{0}+1}+\sum_{i=1}^{k_{0}+1}{k_{0}+1\choose i}\mathbf{z}_{t}^{k_{0}+1-i}r_{2}^{i}(t).
\end{equation}
In a way similar to~\eqref{lm:seq2}, the above implies
\begin{equation}
\label{lm:seq5}
\mathbb{E}\left[\mathbf{z}_{t+1}^{k_{0}+1}\right]\leq (1-\overline{r}_{1}(t))\mathbb{E}\left[\mathbf{z}_{t}^{k_{0}+1}\right]+\sum_{i=1}^{k_{0}+1}{k_{0}+1\choose i}\mathbb{E}\left[\mathbf{z}_{t}^{k_{0}+1-i}\right]r_{2}^{i}(t).
\end{equation}
By the induction hypothesis and the assumptions on the sequence $\{r_{2}(t)\}$, there exist constants $c_{i}$ for $i=1,\cdots,k_{0}+1$, such that,
\begin{equation}
\label{lm:seq6}
\mathbb{E}\left[\mathbf{z}_{t}^{k_{0}+1-i}\right]r_{2}^{i}(t)\leq\frac{c_{i}}{(t+1)^{(k_{0}+1-i)(\delta_{2}-\delta_{1}-\varepsilon_{0})+i\delta_{2}}}=\frac{c_{i}}{(t+1)^{(k_{0}+1)(\delta_{2}-\delta_{1}-\varepsilon_{0})+i(\delta_{1}+\varepsilon_{0})}}
\end{equation}
for all $i=1,\cdots,k_{0}+1$. It is readily seen that the smallest decay rate in the above is attained at $i=1$. Hence, from~\eqref{lm:seq5}-\eqref{lm:seq6}, there exists another constant $c_{0}$, such that,
\begin{equation}
\label{lm:seq7}
\mathbb{E}\left[\mathbf{z}_{t+1}^{k_{0}+1}\right]\leq (1-\overline{r}_{1}(t))\mathbb{E}\left[\mathbf{z}_{t}^{k_{0}+1}\right]+\frac{c_{0}}{(t+1)^{(k_{0}+1)(\delta_{2}-\delta_{1}-\varepsilon_{0})+(\delta_{1}+\varepsilon_{0})}}.
\end{equation}
The deterministic sequence $\left\{\mathbb{E}\left[\mathbf{z}_{t}^{k_{0}+1}\right]\right\}$ then falls under the purview of Lemma~\ref{lm:JSTSP-det} (by taking $\delta_{2}\doteq (k_{0}+1)(\delta_{2}-\delta_{1}-\varepsilon_{0})+(\delta_{1}+\varepsilon_{0})$ and $\delta_{1}\doteq\delta_{1}$. Since $\varepsilon_{0}>0$, an immediate application of Lemma~\ref{lm:JSTSP-det} gives
\begin{equation}
\label{lm:seq8}
\lim_{t\rightarrow\infty}(t+1)^{(k_{0}+1)(\delta_{2}-\delta_{1}-\varepsilon_{0})}\mathbb{E}\left[\mathbf{z}_{t}^{k_{0}+1}\right]=0
\end{equation}
and the induction step follows. This establishes the desired claim in~\eqref{lm:seq1}.

We now complete the proof of Lemma~\ref{lm:seq-gen}. To this end, choose $\overline{\delta}$, such that $0<\overline{\delta}<\delta_{2}-\delta_{1}-\delta_{0}$.
Let $k_{\delta_{0}}$ be a positive integer, such that $k_{\delta_{0}}(\delta_{2}-\delta_{1}-\delta_{0}-\overline{\delta})>1$. Then, for every $\varepsilon>0$, we have
\begin{equation}
\label{lm:seq10}
\mathbb{P}\left((t+1)^{\delta_{0}}\mathbf{z}_{t}>\varepsilon\right)\leq\frac{\mathbb{E}[\mathbf{z}_{t}^{k_{\delta_{0}}}]}{\varepsilon^{k_{\delta_{0}}}(t+1)^{-k_{\delta_{0}}\delta_{0}}}\leq \frac{c}{\varepsilon^{k_{\delta_{0}}}(t+1)^{k_{\delta_{0}}(\delta_{2}-\delta_{1}-\delta_{0}-\overline{\delta})}}.
\end{equation}
The last step is a consequence of the claim in~\eqref{lm:seq1}, by which there exists a constant $c>0$, such that,
\begin{equation}
\label{lm:seq11}
\mathbb{E}[\mathbf{z}_{t}^{k_{\delta_{0}}}]\leq\frac{c}{(t+1)^{k_{\delta_{0}}(\delta_{2}-\delta_{1}-\overline{\delta})}}
\end{equation}
for all $t\geq 0$. Since $k_{\delta_{0}}(\delta_{2}-\delta_{1}-\delta_{0}-\overline{\delta})>1$ by choice, the rightmost term in~\eqref{lm:seq10} is summable in $t$. We thus obtain $\sum_{t=0}^{\infty}\mathbb{P}\left((t+1)^{\delta_{0}}\mathbf{z}_{t}>\varepsilon\right)<\infty$, and, hence,
\begin{equation}
\label{lm:seq13}
\mathbb{P}\left((t+1)^{\delta_{0}}\mathbf{z}_{t}>\varepsilon~~\mbox{i.o.}\right)=0
\end{equation}
by the Borel-Cantelli lemma (i.o. stands for infinitely often in~\eqref{lm:seq13}). Since~\eqref{lm:seq13} holds for arbitrary $\varepsilon>0$, we conclude
that $(t+1)^{\delta_{0}}\mathbf{z}_{t}\rightarrow 0$ a.s. as $t\rightarrow\infty$.
\end{proof}

\begin{proof}[Proof of Lemma~\ref{lm:mean-conv}] Fix $\delta\in\left(0,\delta_{2}-\delta_{1}-\delta_{0}-\frac{1}{2+\Vap_{1}}\right)$. The following is readily verified:

For every $\Vap_{3}>0$, there exists $R_{\Vap_{3}}>0$, such that
\begin{equation}
\label{lm:mean-conv7}
\mathbb{P}\left(\sup_{t\geq 0}\frac{1}{(t+1)^{\frac{1}{2+\Vap_{1}}+\delta}}\left\|U_{t}(1+J_{t})\right\|<R_{\Vap_{3}}\right)>1-\Vap_{3}.
\end{equation}
Indeed, for any $\Vap_{2}>0$, we note that
\begin{eqnarray}
\label{lm:mean-conv8}
\mathbb{P}\left(\frac{1}{(t+1)^{\frac{1}{2+\Vap_{1}}+\delta}}\|J_{t}\|>\Vap_{2}\right) & \leq & \frac{1}{\Vap_{2}^{2+\Vap_{1}}(t+1)^{1+\delta(2+\Vap_{1})}}\mathbb{E}\left[\|J_{t}\|^{2+\Vap_{1}}\right]\nonumber \\ & \leq & \frac{\kappa}{\Vap_{2}^{2+\Vap_{1}}(t+1)^{1+\delta(2+\Vap_{1})}}.
\end{eqnarray}
Since $\delta>0$, the term on the right hand side of~\eqref{lm:mean-conv8} is summable, and by the Borel-Cantelli lemma we may conclude that
\begin{equation}
\label{lm:mean-conv9}
\mathbb{P}\left(\frac{1}{(t+1)^{\frac{1}{2+\Vap_{1}}+\delta}}\|J_{t}\|>\Vap_{2}~\mbox{i.o.}\right)=0.
\end{equation}
Since $\Vap_{2}$ is arbitrary, it follows that
\begin{equation}
\label{lm:mean-conv10}
\mathbb{P}\left(\lim_{t\rightarrow\infty}\frac{1}{(t+1)^{\frac{1}{2+\Vap_{1}}+\delta}}\|J_{t}\|=0\right)=1.
\end{equation}
From the boundedness of $\{U_{t}\}$ and~\eqref{lm:mean-conv10} we may further conclude that
\begin{equation}
\label{lm:mean-conv11}
\mathbb{P}\left(\lim_{t\rightarrow\infty}\frac{1}{(t+1)^{\frac{1}{2+\Vap_{1}}+\delta}}\left\|U_{t}(1+J_{t})\right\|=0\right)=1.
\end{equation}
By Egorov's theorem the a.s. convergence in~\eqref{lm:mean-conv11} is uniform except on a set of arbitrarily small measure, which verifies the claim in~\eqref{lm:mean-conv7}.

We now establish the desired result by a truncation argument. For a
scalar $a$, define its truncation $(a)_{C}$ at level $C>0$ by
\begin{equation}
\label{lm:mean-conv12} (a)_{C}=\left\{\begin{array}{ll}
                                \frac{a}{|a|}\min(|a|,C) & \mbox{if
                                $a\neq 0$}\\
                                0 & \mbox{if $a=0$.}
                                \end{array}
                                \right.
\end{equation}
For a vector, the truncation operation applies component-wise. Now, for each $C>0$, consider the sequence $\{\widehat{\mathbf{z}}_{C}(t)\}$ given by the recursion
\begin{equation}
\label{lm:mean-conv13}
\widehat{\mathbf{z}}_{C}(t+1)=(1-r_{1}(t))\widehat{\mathbf{z}}_{C}(t)+r_{2}(t)\left(U_{t}(1+J_{t})\right)_{C(t+1)^{\frac{1}{2+\Vap_{1}}+\delta}}
\end{equation}
with $\widehat{\mathbf{z}}_{C}(0)=\mathbf{z}_{0}$. Using~\eqref{lm:mean-conv12}, we have
\begin{equation}
\label{lm:mean-conv14}
\widehat{\mathbf{z}}_{C}(t+1)\leq (1-r_{1}(t))\widehat{\mathbf{z}}_{C}(t)+\widehat{r}_{2}(t),
\end{equation}
where
\begin{equation}
\label{lm:mean-conv15}
\widehat{r}_{2}(t)\leq\frac{k_{1}}{(t+1)^{\delta_{2}-\delta-\frac{1}{2+\Vap_{1}}}},~~\forall t
\end{equation}
for some constant $k_{1}>0$. By construction the process $\{\widehat{\mathbf{z}}_{C}(t)\}$ is $\{\mathcal{F}_{t}\}$ adapted and, hence, the recursion in~\eqref{lm:mean-conv14}-\eqref{lm:mean-conv15} falls under the purview of Lemma~\ref{lm:seq-gen}. Thus, for every $C>0$, we have $(t+1)^{\delta_{0}}\widehat{\mathbf{z}}_{C}(t)\rightarrow 0$ a.s. as $t\rightarrow\infty$, since $\delta_{0}<\delta_{2}-\delta_{1}-\delta-\frac{1}{2+\Vap_{1}}$.

Now, for $\Vap_{3}>0$, consider the sequence $\{\widehat{\mathbf{z}}_{R_{\Vap_{3}}}(t)\}$, where $R_{\Vap_{3}}>0$ is the constant in~\eqref{lm:mean-conv7}. Using~\eqref{lm:mean-conv7} and~\eqref{lm:mean-conv13} we may conclude that
\begin{equation}
\label{lm:mean-conv17}
\mathbb{P}\left(\inf_{t\geq 0}\left(\widehat{\mathbf{z}}_{R_{\Vap_{3}}}(t)-\mathbf{z}_{t}\right)\geq 0\right)>1-\Vap_{3}.
\end{equation}
Since all processes involved are non-negative, it readily follows from~\eqref{lm:mean-conv17} that
\begin{equation}
\label{lm:mean-conv18}
\mathbb{P}\left(\lim_{t\rightarrow\infty}(t+1)^{\delta_{0}}\mathbf{z}_{t}=0\right)>1-\Vap_{3}.
\end{equation}
The lemma follows by taking $\Vap_{3}$ to zero in~\eqref{lm:mean-conv18}.
\end{proof}

\begin{proof}[Proof of Lemma~\ref{lm:conn}] Let $\mathcal{L}$ denote the set of possible Laplacian matrices (necessarily finite) and $\mathcal{D}$ the distribution on $\mathcal{L}$ induced by the link formation process. Since the set of Laplacian matrices is finite, the set $\mathcal{L}$ may be chosen such that $\underline{p}=\inf_{L\in\mathcal{L}}p_{L}>0$, with $p_{L}=\mathbb{P}(L_{t}=L)$ for each $L\in\mathcal{L}$ and $\sum_{L\in\mathcal{L}}p_{L}=1$. The hypothesis $\lambda_{2}(\overline{L})>0$ implies that for every $\mathbf{z}\in\mathcal{C}^{\perp}$,
\begin{equation}
\label{lm:conn7}
\sum_{L\in\mathcal{L}}\mathbf{z}^{T}L\mathbf{z}\geq\sum_{L\in\mathcal{L}}\mathbf{z}^{T}(p_{L}L)\mathbf{z}=\mathbf{z}^{T}\overline{L}\mathbf{z}\geq\lambda_{2}(\OL)\|\mathbf{z}\|^{2}.
\end{equation}
Denoting by $|\mathcal{L}|$ the cardinality of $\mathcal{L}$, it follows from~\eqref{lm:conn7} that for each $\mathbf{z}\in\mathcal{C}^{\perp}$ there exists some $L_{\mathbf{z}}\in\mathcal{L}$, such that $\mathbf{z}^{T}L_{\mathbf{z}}\mathbf{z}\geq (\lambda_{2}(\OL)/|\mathcal{L}|)\|\mathbf{z}\|^{2}$. Moreover, since the set $\mathcal{L}$ is finite, the mapping $L_{\mathbf{z}}:\mathcal{C}^{\perp}\longrightarrow\mathcal{L}$ may be realized as a measurable function.

For each $L\in\mathcal{L}$, the eigenvalues of the matrix $I_{NM}-\beta_{t}L\otimes I_{M}$ are $1$ and $1-\beta_{t}\lambda_{n}(L)$, $2\leq n\leq N$, each being repeated $M$ times. Hence, for $t\geq t_{0}$ (large enough), $\|I_{NM}-\beta_{t}L\otimes I_{M}\|\leq 1$ and $\|(I_{NM}-\beta_{t}L\otimes I_{M})\mathbf{z}\|\leq\|\mathbf{z}\|$ for every $\mathbf{z}\in\mathbb{R}^{NM}$. Hence, the functional $r_{L,\mathbf{z}}$  given by
\begin{equation}
\label{lm:conn8}r_{L,\mathbf{z}}=\left\{\begin{array}{ll}
                                1 & \mbox{if $t<t_{0}$ or $\mathbf{z}=\mathbf{0}$}\\
                                1-\frac{\|(I_{NM}-\beta_{t}L\otimes I_{M})\mathbf{z}\|}{\|\mathbf{z}\|} & \mbox{otherwise}
                                \end{array}
                                \right.
\end{equation}
is jointly measurable in $L$ and $\mathbf{z}$ and
satisfies $0\leq r_{L,\mathbf{z}}\leq 1$ for each pair $(L,\mathbf{z})$. Let $\{r_{t}\}$ be the $\{\mathcal{F}_{t+1}\}$ adapted process given by $r_{t}=r_{L_{t},\mathbf{z}_{t}}$ for each $t$, and $\|(I_{NM}-\beta_{t}L\otimes I_{M})\mathbf{z}_{t}\|=(1-r_{t})\|\mathbf{z}_{t}\|$ a.s. for each $t$. We now need to verify that $\{r_{t}\}$ satisfies~\eqref{lm:conn2} for some $c_{r}>0$. To this end,  for $t$ large enough,
\begin{eqnarray}
\label{lm:conn9}\|(I_{NM}-\beta_{t}L_{\mathbf{z}_{t}}\otimes I_{M})\mathbf{z}_{t}\|^{2} & = & \mathbf{z}_{t}^{T}(I_{NM}-2\beta_{t}L_{\mathbf{z}_{t}}\otimes I_{M})\mathbf{z}_{t}+\beta_{t}^{2}\mathbf{z}_{t}^{T}(L_{\mathbf{z}_{t}}\otimes I_{M})^{2}\mathbf{z}_{t}\nonumber \\ & \leq & \left(1-2\beta_{t}\lambda_{2}(\OL)/|\mathcal{L}|\right)\|\mathbf{z}_{t}\|^{2}+c_{1}\beta_{t}^{2}\|\mathbf{z}_{t}\|^{2}\nonumber \\ & \leq & \left(1-\beta_{t}\lambda_{2}(\OL)/|\mathcal{L}|\right)\|\mathbf{z}_{t}\|^{2},
\end{eqnarray}
where we have used the definition of the function $L_{\mathbf{z}}$, the boundedness of the Laplacian matrix and the fact that $\beta_{t}\rightarrow 0$. Hence, by making $t_{0}$ larger if necessary, we have
\begin{equation}
\label{lm:conn10}
\|(I_{NM}-\beta_{t}L_{\mathbf{z}_{t}}\otimes I_{M})\mathbf{z}_{t}\|\leq\left(1-\beta_{t}\frac{\lambda_{2}(\OL)}{4|\mathcal{L}|}\right)\|\mathbf{z}_{t}\|
\end{equation}
for all $t\geq t_{0}$. Now, by~\eqref{lm:conn10}
\begin{eqnarray}
\label{lm:conn11}
\mathbb{E}\left[\|(I_{NM}-\beta_{t}L\otimes I_{M})\mathbf{z}_{t}\|~|~\mathcal{F}_{t}\right] & = & \sum_{L\in\mathcal{L}}p_{L}\left(1-r_{L,\mathbf{z}_{t}}\right)\|\mathbf{z}_{t}\|\nonumber \\ & \leq & \left(1-\left(\underline{p}\beta_{t}\frac{\lambda_{2}(\OL)}{4|\mathcal{L}|}+\sum_{L\neq L_{\mathbf{z}_{t}}}p_{L}r_{L,\mathbf{z}_{t}}\right)\right)\|\mathbf{z}_{t}\|.
\end{eqnarray}
Since $\sum_{L\neq L_{\mathbf{z}_{t}}}p_{L}r_{L,\mathbf{z}_{t}}\geq 0$, we have for $t\geq t_{0}$,
\begin{equation}
\label{lm:conn12}
(1-\mathbb{E}[r_{t}|\mathcal{F}_{t}])\|\mathbf{z}_{t}\|=\mathbb{E}\left[\|(I_{NM}-\beta_{t}L\otimes I_{M})\mathbf{z}_{t}\|~|~\mathcal{F}_{t}\right]\leq \left(1-\underline{p}\beta_{t}\frac{\lambda_{2}(\OL)}{4|\mathcal{L}|}\right)\|\mathbf{z}_{t}\|.
\end{equation}
Since, by definition $r_{t}=1$ on the set $\{\mathbf{z}_{t}=\mathbf{0}\}$, it follows that
\begin{equation}
\label{lm:conn13}
\mathbb{E}[r_{t}|\mathcal{F}_{t}]\geq \frac{\underline{p}\lambda_{2}(\OL)}{4|\mathcal{L}|}\beta_{t}
\end{equation}
for all $t\geq t_{0}$, thus establishing the assertion.
\end{proof}

}

\section{Proofs in Section~\ref{sec:conv-asym}}
\label{app:est}
{\small
\begin{proof}[Proof of Lemma~\ref{lm:Gcons}] We will show the desired convergence in the matrix Frobenius norm (denoted by $\|\cdot\|_{F}$ in the following). Since the matrix space under consideration is finite dimensional, the convergence in $\mathcal{L}_{2}$ norm will follow. The existence of quadratic moments implies the convergence of the sample covariances (see~\eqref{Q_up}) to the true covariances and, hence, for each $n$, $Q_{n}(t)\rightarrow R_{n}$ a.s. Since, in addition, the sequence $\{\gamma_{t}\}$ in~\eqref{gain_up} goes to zero, we may choose an a.s. finite random variable $R_{2}$, such that for each $n$,
\begin{equation}
\label{lm:Gb2}
\Past\left(\sup_{t\geq 0}\left\|H_{n}^{T}\left(Q_{n}(t)+\gamma_{t}I_{M_{n}}\right)^{-1}H_{n}\right\|\leq R_{2}<\infty\right)=1.
\end{equation}
By construction, the matrix sequences $\{G_{n}(t)\}$ and $\{Q_{n}(t)\}$ are symmetric for each $n$. Let $\widetilde{G}_{n}(t)=G_{n}(t)-\Gavg(t)$ denote the deviation of the Grammian estimate at agent $n$ from the instantaneous network average $\Gavg(t)$. Also, let $\widetilde{G}_{t}$ and $D_{t}$ respectively denote the matrices $[\widetilde{G}_{1}(t),\cdots,\widetilde{G}_{N}(t)]^{T}$ and $[D_{1}(t),\cdots,D_{N}(t)]^{T}$, where $D_{n}(t)=\left(Q_{n}(t)+\gamma_{t}I_{M_{n}}\right)^{-1}$ for each $n$. Using the following readily verifiable properties of the Laplacian:
\begin{equation}
\label{lm:Gcons3}
\left(\mathbf{1}_{N}\otimes I_{M}\right)^{T}\left(L_{t}\otimes I_{M}\right)=0,~~~\left(L_{t}\otimes I_{M}\right)\left(\mathbf{1}_{N}\otimes\Gavg(t)\right)=\mathbf{0},
\end{equation}
we have
\begin{equation}
\label{lm:Gcons5}
\widetilde{G}_{t+1}=\left(I_{NM}-\beta_{t}\left(L_{t}\otimes I_{M}\right)-\alpha_{t}I_{NM}\right)\widetilde{G}_{t}+\alpha_{t}\left(\left(D_{t}-\Davg(t)\right)\right),
\end{equation}
where $\Davg(t)=\frac{1}{N}\sum_{n=1}^{N}D_{n}(t)$. Note that, by~\eqref{lm:Gb2}, there exists an $\{\mathcal{F}_{t}\}$ adapted a.s. bounded process $\{U_{t}\}$, such that $\sup_{t\geq 0}\|D_{t}-\Davg(t)\|_{F}\leq U_{t}$ a.s. For $m\in\{1,\cdots,M\}$, let $\widetilde{G}_{m,t}$ denote the $m$-th column of $\widetilde{G}_{t}$. The process $\{\widetilde{G}_{m,t}\}$ is $\{\mathcal{F}_{t}\}$ adapted and $\widetilde{G}_{m,t}\in\mathcal{C}^{\perp}$ for each $t$. Then, by Lemma~\ref{lm:conn} there exists a $[0,1]$-valued $\{\mathcal{F}_{t+1}\}$ adapted process $\{r_{m,t}\}$, such that,
\begin{equation}
\label{lm:Gcons200}
\|(I_{NM}-\beta_{t}L_{t}\otimes I_{M})\widetilde{G}_{m,t}\|\leq (1-r_{m,t})\|\widetilde{G}_{m,t}\|
\end{equation}
and $\East[r_{m,t}|\mathcal{F}_{t}]\geq c_{m,r}/(t+1)^{\tau_{2}}$ a.s. for $t\geq t_{0}$ sufficiently large. Noting that the square of the Frobenius norm is the sum of the squared column $\mathcal{L}_{2}$ norms, we have
\begin{equation}
\label{lm:Gcons201}
\|(I_{NM}-\beta_{t}L_{t}\otimes I_{M})\widetilde{G}_{t}\|_{F}^{2}\leq\sum_{m=1}^{M}(1-r_{m,t})^{2}\|\widetilde{G}_{m,t}\|^{2}\leq (1-r_{t})^{2}\|\widetilde{G}_{t}\|^{2}_{F},
\end{equation}
where $\{r_{t}\}$ is the $\{\mathcal{F}_{t+1}\}$ adapted process given by $r_{t}=r_{1,t}\wedge r_{2,t}\wedge\cdots\wedge r_{M,t}$. By the conditional Jensen's inequality, we obtain
\begin{equation}
\label{lm:conn202}
\East[r_{t}|\mathcal{F}_{t}]\geq\wedge_{m=1}^{M}\East[r_{m,t}|\mathcal{F}_{t}]\geq c_{r}/(t+1)^{\tau_{2}}
\end{equation}
for some $c_{r}>0$ and $t\geq t_{0}$. Recall $\{\alpha_{t}\}$ from~\eqref{weight}. Using~\eqref{lm:Gcons201}, we finally get
\begin{align}
\nonumber
\|(I_{NM}-\beta_{t}L_{t}\otimes I_{M}-\alpha_{t}I_{NM})\widetilde{G}_{t}\|_{F}
 \leq &
\|(I_{NM}-\beta_{t}L_{t}\otimes I_{M})\widetilde{G}_{t}\|_{F}+\alpha_{t}\|\widetilde{G}_{t}\|_{F}
\\
 \leq & (1-r_{t})\|\widetilde{G}_{t}\|_{F}+\alpha_{t}\|\widetilde{G}_{t}\|_{F} \nonumber
 \\
 \leq &
\left(1-r_{t}/2\right)\|\widetilde{G}_{t}\|_{F}
\label{lm:Gcons203}
\end{align}
for $t\geq t_{0}$. From~\eqref{lm:Gcons5} and~\eqref{lm:Gcons203} we then have
\begin{equation}
\label{lm:Gcons7}
\|\widetilde{G}_{t+1}\|_{F}\le\|(I_{NM}-\beta_{t}L_{t}\otimes I_{M}-\alpha_{t}I_{NM})\widetilde{G}_{t}\|_{F}+\alpha_{t}U_{t}\leq\left(1-r_{t}/2\right)\|\widetilde{G}_{t}\|_{F}+\alpha_{t}U_{t}.
\end{equation}
By~\eqref{lm:conn202} and since $\beta_{t}/\alpha_{t}\rightarrow\infty$ as $t\rightarrow\infty$, the recursion in~\eqref{lm:Gcons7} clearly falls under the purview of Lemma~\ref{lm:mean-conv}, and we conclude that
$\|\widetilde{G}_{t}\|_{F}\rightarrow 0$ a.s. as $t\rightarrow\infty$. The convergence in the $\mathcal{L}_{2}$ norm follows immediately.
\end{proof}
\begin{proof}[Proof of Lemma~\ref{lm:Gavg}]
The process $\{\Gavg(t)\}$ satisfies the following recursion:
\begin{equation}
\label{lm:Gavg2}
\Gavg(t+1)=(1-\alpha_{t})\Gavg(t)+\alpha_{t}\Davg(t).
\end{equation}
Let $\TGavg(t)$ denote the residual $\Gavg(t)-\OverlineSigmac$ and the process $\{\TGavg(t)\}$ satisfies
\begin{equation}
\label{lm:Gavg3}
\TGavg(t+1)=\left(1-\alpha_{t}\right)\TGavg(t)+\alpha_{t}\left(\Davg(t)-\OverlineSigmac\right).
\end{equation}
By Eqn.~139, Lemma~25 in~\cite{KarMouraRamanan-Est} there exist $t_{0}$ sufficiently large and a constant $B$ such that
\begin{equation}
\label{lm:Gavg4}
0\leq \sum_{k=s}^{t-1}\left((
\prod_{l=k+1}^{t-1}(1-\alpha_{l}))\alpha_{k}\right)\leq
B,
\end{equation}
for all positive integers $t$ and $s$ with $t_{0}\leq s\leq t$. Also, the convergence of the sample covariances and the fact that $\gamma_{t}\rightarrow 0$ as $t\rightarrow\infty$ imply $\Davg(T)\rightarrow\OverlineSigmac$ a.s. as $t\rightarrow\infty$. Hence, for a given $\varepsilon>0$, we may choose $t_{\varepsilon}>t_{0}$ such that $\left\|\Davg(t)-\OverlineSigmac\right\|<\varepsilon$ for all $t\geq t_{\varepsilon}$. From~\eqref{lm:Gavg3}, we then have for $t>t_{\varepsilon}$
\begin{align}
\label{lm:Gavg6}
\|\TGavg(t)\|  \leq & \left|(\prod_{k=t_{\varepsilon}}^{t-1}(1-\alpha_{k}))\right|
\left\|\TGavg(t_{\varepsilon})\right\|+\sum_{k=t_{\varepsilon}}^{t-1}\left(\left(
\prod_{l=k+1}^{t-1}(1-\alpha_{l})\right)\alpha_{k}\varepsilon\right)\nonumber \\ \leq & \left|(\prod_{k=t_{\varepsilon}}^{t-1}(1-\alpha_{k}))\right|
\left\|\TGavg(t_{\varepsilon})\right\|+B\varepsilon.
\end{align}
Since $\sum_{t\geq 0}\alpha_{t}=\infty$ the first term on the right hand side of~\eqref{lm:Gavg6} goes to zero as $t\rightarrow\infty$, and we have
$\limsup_{t\rightarrow\infty}\|\TGavg(t)\|\leq B\varepsilon$.
Since $\varepsilon>0$ is arbitrary, we conclude that $\TGavg(t)\rightarrow\mathbf{0}$ a.s. as $t\rightarrow\infty$ by taking $\varepsilon$ to zero. The desired assertion follows immediately.
\end{proof}

\begin{proof}[Proof of Proposition~\ref{prop:est}] A version of this result was established in~\cite{JSTSP-Kar-Moura} (Lemma 6) for the case of constant gains $K_{n}(t)$. In the following we generalize the arguments of~\cite{JSTSP-Kar-Moura} to time-varying adaptive gains. To this end we show
\begin{equation}
\label{prop:est2}
\inf_{\|\mathbf{z}\|=1}\mathbf{z}^{T}\left(\frac{\beta_{t}}{\alpha_{t}}\overline{L}\otimes I_{M}+\mathcal{K}\mathcal{H}\right)\mathbf{z}>0
\end{equation}
for all $t$ sufficiently large, where $\mathcal{K}=\diag\left(K_{1},\cdots,K_{N}\right)$.

A vector $\mathbf{z}\in\mathbb{R}^{NM}$ may be decomposed as $\mathbf{z}=\mathbf{z}_{\C}+\mathbf{z}_{\PC}$, with $\mathbf{z}_{\C}$ denoting its projection on the consensus or agreement subspace $\mathcal{C}$,
\begin{equation}
\label{prop:est3}
\mathcal{C}=\left\{\mathbf{z}\in\mathbb{R}^{NM}~|~\mathbf{z}=\mathbf{1}_{N}\otimes\mathbf{a}~\mbox{for some $\mathbf{a}\in\mathbb{R}^{M}$}\right\},
\end{equation}
and $\mathbf{z}_{\PC}$ the orthogonal complement. Also, denoting by $\mathcal{D}_{\mathcal{K}}$ the symmetricized version of $\mathcal{K}\mathcal{H}$, i.e.,
$\mathcal{D}_{\mathcal{K}}=\frac{1}{2}\left(\mathcal{K}\mathcal{H}+\mathcal{H}^{T}\mathcal{K}^{T}\right)$,
standard matrix manipulations and properties of the Laplacian yield
\begin{equation}
\label{prop:est5}
\mathbf{z}^{T}\left(\frac{\beta_{t}}{\alpha_{t}}\overline{L}\otimes I_{M}+\mathcal{K}\mathcal{H}\right)\mathbf{z}\geq\frac{\beta_{t}}{\alpha_{t}}\lambda_{2}(\overline{L})\left\|\mathbf{z}_{\PC}\right\|^{2}+\mathbf{z}_{\PC}^{T}\mathcal{D}_{\mathcal{K}}\mathbf{z}_{\PC}+2\mathbf{z}_{\C}^{T}\mathcal{D}_{\mathcal{K}}\mathbf{z}_{\PC}+\mathbf{z}_{\C}^{T}\mathcal{D}_{\mathcal{K}}\mathbf{z}_{\C}.
\end{equation}
By construction, $\sum_{n=1}^{N}K_{n}H_{n}=\OverlineSigmac^{-1}\sum_{n=1}^{N}H_{n}^{T}R_{n}^{-1}H_{n}=NI_{M}$, and, hence, we note that $\mathbf{z}_{\C}^{T}\mathcal{D}_{\mathcal{K}}\mathbf{z}_{\C} = \left\|\mathbf{z}_{\C}\right\|^{2}$ for each $\mathbf{z}\in\mathbb{R}^{NM}$.
Let us choose a constant $c_{1}>0$ such that
\begin{equation}
\label{prop:est8}
\mathbf{z}_{\PC}^{T}\mathcal{D}_{\mathcal{K}}\mathbf{z}_{\PC}\geq -c_{1}\left\|\mathbf{z}_{\PC}\right\|^{2}~~\mbox{and}~~\mathbf{z}_{\C}^{T}\mathcal{D}_{\mathcal{K}}\mathbf{z}_{\PC}\geq -c_{1}\left\|\mathbf{z}_{\C}\right\|\left\|\mathbf{z}_{\PC}\right\|.
\end{equation}
It then follows from~\eqref{prop:est5} that
\begin{equation}
\label{prop:est10}
\mathbf{z}^{T}\left(\frac{\beta_{t}}{\alpha_{t}}\overline{L}\otimes I_{M}+\mathcal{K}\mathcal{H}\right)\mathbf{z}\geq\left(\frac{\beta_{t}}{\alpha_{t}}\lambda_{2}(\overline{L})-c_{1}\right)\left\|\mathbf{z}_{\PC}\right\|^{2}-2c_{1}\left\|\mathbf{z}_{\C}\right\|\left\|\mathbf{z}_{\PC}\right\|+\left\|\mathbf{z}_{\C}\right\|^{2}.
\end{equation}
Since $\beta_{t}/\alpha_{t}\rightarrow\infty$ and $\lambda_{2}(\overline{L})>0$, there exists $t_{1}$ sufficiently large such that
\begin{equation}
\label{prop:est11}
\frac{\beta_{t}}{\alpha_{t}}\lambda_{2}(\overline{L})-c_{1}>c_{1}^{2},~~~\forall t\geq t_{1}.
\end{equation}
We now verify~\eqref{prop:est2} for $t\geq t_{1}$. To this end, assume $\|\mathbf{z}\|=1$. In case $\mathbf{z}_{\C}=\mathbf{0}$ ($\|\mathbf{z}_{\PC}\|=1$), we have from~\eqref{prop:est10}
\begin{equation}
\label{prop:est12}
\mathbf{z}^{T}\left(\frac{\beta_{t}}{\alpha_{t}}\overline{L}\otimes I_{M}+\mathcal{K}\mathcal{H}\right)\mathbf{z}\geq\frac{\beta_{t}}{\alpha_{t}}\lambda_{2}(\overline{L})-c_{1}>0.
\end{equation}
For the other case, i.e., $\mathbf{z}_{\C}\neq\mathbf{0}$,
\begin{equation}
\label{prop:est13}
\mathbf{z}^{T}\left(\frac{\beta_{t}}{\alpha_{t}}\overline{L}\otimes I_{M}+\mathcal{K}\mathcal{H}\right)\mathbf{z}\geq\left\|\mathbf{z}_{\C}\right\|^{2}\left[\left(\frac{\beta_{t}}{\alpha_{t}}\lambda_{2}(\overline{L})-c_{1}\right)\frac{\|\mathbf{z}_{\PC}\|^{2}}{\|\mathbf{z}_{\C}\|^{2}}-2c_{1}\frac{\|\mathbf{z}_{\PC}\|}{\|\mathbf{z}_{\C}\|}+1\right]>0,
\end{equation}
where the last inequality follows from the fact that the quadratic functional of $\frac{\|\mathbf{z}_{\PC}\|}{\|\mathbf{z}_{\C}\|}$ is always positive due to the discriminant condition imposed by~\eqref{prop:est11}. We thus conclude that
\begin{equation}
\label{prop:est14}
\mathbf{z}^{T}\left(\frac{\beta_{t}}{\alpha_{t}}\overline{L}\otimes I_{M}+\mathcal{K}\mathcal{H}\right)\mathbf{z}>0
\end{equation}
for all $t\geq t_{1}$ and $\mathbf{z}$, such that $\|\mathbf{z}\|=1$. Since the quadratic form in~\eqref{prop:est14} is a continuous function on the compact unit circle, we may further conclude that
\begin{equation}
\label{prop:est16}
\inf_{\|\mathbf{z}\|=1}\mathbf{z}^{T}\left(\frac{\beta_{t}}{\alpha_{t}}\overline{L}\otimes I_{M}+\mathcal{K}\mathcal{H}\right)\mathbf{z}>c_{2}>0,
\end{equation}
for some positive constant $c_{2}$, thus verifying the assertion in~\eqref{prop:est2} for all $t\geq t_{1}$. To complete the proof of Proposition~\ref{prop:est}, choose any $0<\Vap<c_{2}$. It then follows from~\eqref{prop:est16} and continuity that for $t\geq t_{1}$ and arbitrary $\mathbf{z}\in\mathbb{R}^{NM}$,
\begin{eqnarray}
\label{prop:est18}
\mathbf{z}^{T}\left(\beta_{t}\overline{L}\otimes I_{M}+\alpha_{t}\widetilde{\mathcal{K}}\mathcal{H}\right)\mathbf{z} & \geq &\alpha_{t}\left\|\mathbf{z}\right\|^{2}\left[\inf_{\|\overline{\mathbf{z}}\|=1}\overline{\mathbf{z}}^{T}\left(\frac{\beta_{t}}{\alpha_{t}}\overline{L}\otimes I_{M}+\widetilde{\mathcal{K}}\mathcal{H}\right)\overline{\mathbf{z}}\right]\nonumber \\ & \geq & \left(c_{2}-\Vap\right)\alpha_{t}\left\|\mathbf{z}\right\|^{2},
\end{eqnarray}
thus verifying the assertion of Proposition~\ref{prop:est} with $\Vap_{\mathcal{K}}=\Vap$, $t_{\mathcal{K}}=t_{1}$ and $c_{\mathcal{K}}=c_{2}-\Vap$.
\end{proof}

\begin{proof}[Proof of Proposition~\ref{prop:cons}] By~\eqref{prop:est10} in Proposition~\ref{prop:est}, there exists a constant $c_{1}>0$ such that for arbitrary $\mathbf{z}\in\mathbb{R}^{NM}$
\begin{equation}
\label{prop:cons3}
\mathbf{z}^{T}\left(\frac{\beta_{t}}{\alpha_{t}}\overline{L}\otimes I_{M}+\mathcal{K}\mathcal{H}\right)\mathbf{z}\geq\left(\frac{\beta_{t}}{\alpha_{t}}\lambda_{2}(\overline{L})-c_{1}\right)\left\|\mathbf{z}_{\PC}\right\|^{2}-2c_{1}\left\|\mathbf{z}_{\C}\right\|\left\|\mathbf{z}_{\PC}\right\|+\left\|\mathbf{z}_{\C}\right\|^{2}.
\end{equation}
Hence for $\widetilde{\mathcal{K}}$ satisfying~\eqref{prop:cons2}, we have
\begin{eqnarray}
\label{prop:cons4}
\mathbf{z}^{T}\left(\frac{\beta_{t}}{\alpha_{t}}\overline{L}\otimes I_{M}+\widetilde{\mathcal{K}}\mathcal{H}\right)\mathbf{z} & \geq & \mathbf{z}^{T}\left(\frac{\beta_{t}}{\alpha_{t}}\overline{L}\otimes I_{M}+\mathcal{K}\mathcal{H}\right)\mathbf{z}-\Vap\left\|\mathbf{z}\right\|^{2}\nonumber \\ & = & \left(\frac{\beta_{t}}{\alpha_{t}}\lambda_{2}(\overline{L})-c_{1}-\Vap\right)\left\|\mathbf{z}_{\PC}\right\|^{2}-2c_{1}\left\|\mathbf{z}_{\C}\right\|\left\|\mathbf{z}_{\PC}\right\|+\left(1-\Vap\right)\left\|\mathbf{z}_{\C}\right\|^{2}.
\end{eqnarray}
Using the fact that $0<\Vap<1$, we have
\begin{eqnarray}
\label{prop:cons5}
\mathbf{z}^{T}\left(\frac{\beta_{t}}{\alpha_{t}}\overline{L}\otimes I_{M}+\widetilde{\mathcal{K}}\mathcal{H}\right)\mathbf{z} & \geq & \left(\frac{\beta_{t}}{2\alpha_{t}}\lambda_{2}(\overline{L})+\left(\frac{\beta_{t}}{2\alpha_{t}}\lambda_{2}(\overline{L})-c_{1}-\Vap-\frac{c_{1}^{2}}{1-\Vap}\right)\right)\left\|\mathbf{z}_{\PC}\right\|^{2}\nonumber \\ & & +\left(\frac{c_{1}}{\sqrt{1-\Vap}}\left\|\mathbf{z}_{\PC}\right\|-\sqrt{1-\Vap}\left\|\mathbf{z}_{\C}\right\|\right)^{2}.
\end{eqnarray}
Since $\lambda_{2}(\overline{L})>0$ and $\beta_{t}/\alpha_{t}\rightarrow\infty$ as $t\rightarrow\infty$, there exists $t_{\Vap}$ (large enough), such that,
\begin{equation}
\label{prop:cons6}
\left(\frac{\beta_{t}}{2\alpha_{t}}\lambda_{2}(\overline{L})-c_{1}-\Vap-\frac{c_{1}^{2}}{1-\Vap}\right)\geq 0
\end{equation}
for all $t\geq t_{\Vap}$. We may then conclude from~\eqref{prop:cons5} that
\begin{equation}
\label{prop:cons7}
\mathbf{z}^{T}\left(\frac{\beta_{t}}{\alpha_{t}}\overline{L}\otimes I_{M}+\widetilde{\mathcal{K}}\mathcal{H}\right)\mathbf{z}\geq\frac{\beta_{t}}{2\alpha_{t}}\lambda_{2}(\overline{L})\left\|\mathbf{z}_{\PC}\right\|^{2}
\end{equation}
and, hence
\begin{equation}
\label{prop:cons8}
\mathbf{z}^{T}\left(\beta_{t}\overline{L}\otimes I_{M}+\alpha_{t}\widetilde{\mathcal{K}}\mathcal{H}\right)\mathbf{z}\geq\frac{\lambda_{2}(\overline{L})}{2}\beta_{t}\left\|\mathbf{z}_{\PC}\right\|^{2}
\end{equation}
for all $t\geq t_{\Vap}$, $\mathbf{z}\in\mathbb{R}^{NM}$ and $\widetilde{K}$ satisfying~\eqref{prop:cons2}. This establishes the assertion.
\end{proof}
\begin{proof}[Proof of Lemma~\ref{lm:estc}] Let the residual $\widetilde{\mathbf{x}}_{n}(t)=\mathbf{x}_{n}(t)-\Xavg(t)$. Then arguments along the lines of~\eqref{lm:Gcons3}-\eqref{lm:Gcons5} show that the process $\widetilde{\mathbf{x}}_{t}=[\widetilde{\mathbf{x}}_{1}^{T}(t),\cdots,\widetilde{\mathbf{x}}_{N}^{T}(t)]^{T}$ satisfies the recursion
\begin{equation}
\label{lm:estc3}
\widetilde{\mathbf{x}}_{t+1}=\left(I_{NM}-\beta_{t}L_{t}\otimes I_{M}\right)\widetilde{\mathbf{x}}_{t}+\alpha_{t}\widetilde{\mathbf{z}}_{t},
\end{equation}
where the process $\{\widetilde{\mathbf{z}}_{t}\}$ is defined as
\begin{equation}
\label{lm:estc4}
\widetilde{\mathbf{z}}_{t}=\left(I_{NM}-\frac{1}{N}\mathbf{1}_{N}\otimes\left(\mathbf{1}_{N}\otimes I_{M}\right)^{T}\right)\mathcal{K}_{t}\left(\mathbf{y}_{t}-\mathcal{H}\mathbf{x}_{t}\right).
\end{equation}
Since $\mathcal{K}_{t}\rightarrow\mathcal{K}$ as $t\rightarrow\infty$, the process $\{\mathbf{x}_{t}\}$ is bounded (Lemma~\ref{lm:estb}), and the observation noise $\mathbf{\zeta}_{t}$ satisfies~\textbf{(A.5)}, there exist two $\mathbb{R}_{+}$ valued processes: 1) an $\mathcal{F}_{t}$-adapted $\{U_{t}\}$   satisfying $\sup_{t\geq 0}\|U_{t}\|<\infty$ a.s.; and  (2) an i.i.d. $\{\mathcal{F}_{t+1}\}$ adapted $\{J_{t}\}$ independent of $\mathcal{F}_{t}$ for each $t$ and $\East\left[\left\|J_{t}\right\|^{2+\Vap_{1}}\right]<\infty$, such that
\begin{equation}
\label{lm:estc7}
\left\|\widetilde{\mathbf{z}}_{t}\right\|\leq U_{t}\left(1+J_{t}\right).
\end{equation}
Since $\widetilde{\mathbf{x}}_{t}\in\mathcal{C}^{\perp}$ for all $t$, by Lemma~\ref{lm:conn} there exists an $\{\mathcal{F}_{t+1}\}$ adapted $\mathbb{R}_{+}$ valued process $\{r_{t}\}$ with $0\leq r_{t}\leq 1$ a.s. such that
\begin{equation}
\label{lm:estc8}
\left\|\left(I_{NM}-\beta_{t}L_{t}\otimes I_{M}-\mathcal{P}_{NM}\right)\widetilde{\mathbf{x}}_{t}\right\|\leq (1-r_{t})\left\|\widetilde{\mathbf{x}}_{t}\right\|
\end{equation}
for all $t$ (large enough) and a constant $c_{r}>0$ such that for all $t$
\begin{equation}
\label{lm:estc9}
\East\left[r_{t}~|~\mathcal{F}_{t}\right]\geq \frac{c_{r}}{(t+1)^{\tau_{2}}}~~\mbox{a.s.}
\end{equation}
From the above development we conclude that
\begin{equation}
\label{lm:estc10}
\left\|\widetilde{\mathbf{x}}_{t+1}\right\|\leq \left(1-r_{t}\right)\left\|\widetilde{\mathbf{x}}_{t}\right\|+\alpha_{t}U_{t}\left(1+J_{t}\right)
\end{equation}
for all $t$ (large enough). The recursion~\eqref{lm:estc10} clearly falls under the purview of Lemma~\ref{lm:mean-conv},and we have the assertion
\begin{equation}
\label{lm:estc11}
\Past\left(\lim_{t\rightarrow\infty}(t+1)^{\tau_{0}}\widetilde{\mathbf{x}}_{t}=0\right)=1
\end{equation}
for all $\tau_{0}\in\left[0,\tau_{1}-\tau_{2}-\frac{1}{2+\Vap_{1}}\right)$. This establishes the claim.
\end{proof}
}

\textbf{Proof of Lemma~\ref{lm:avgr}}

{\small
We will use the following approximation result from~\cite{Fabian-1} in the proof.
\begin{proposition}[Lemma 4.3 in~\cite{Fabian-1}]
\label{prop:Fab-1}Let $\{b_{t}\}$ be a scalar sequence satisfying
\begin{equation}
\label{prop:Fab-11}
b_{t+1}\leq\left(1-\frac{c}{t+1}\right)b_{t}+d_{t}(t+1)^{-\tau}
\end{equation}
where $c>\tau$, $\tau>0$, and the sequence $\{d_{t}\}$ is summable. Then $\limsup_{t\rightarrow\infty}(t+1)^{\tau}b_{t}<\infty$.
\end{proposition}

The following generalized convergence criterion of dependent stochastic sequences will also be useful.
\begin{proposition}[Lemma 10 in~\cite{Dub-Freed}]
\label{prop:DF} Let $\{\overline{J}_{t}\}$ be an $\mathbb{R}$ valued $\{\mathcal{F}_{t+1}\}$ adapted process such that $\mathbb{E}[\overline{J}_{t}|\mathcal{F}_{t}]=0$ a.s. for each $t\geq 1$. Then the sum $\sum_{t\geq 0}\overline{J}_{t}$ exists and is finite a.s. on the set where $\sum_{t\geq 0}\mathbb{E}[\overline{J}_{t}^{2}|\mathcal{F}_{t}]$ is finite.
\end{proposition}
\begin{proof}[Proof of Lemma~\ref{lm:avgr}]
For each $\delta>0$ recall the construction in~\eqref{lm:avgc2}-\eqref{lm:avgc8}. Clearly, it suffices by the arguments in Lemma~\ref{lm:avgc} to establish the required convergence rate claim for each of the processes $\{\mathbf{z}_{t}^{\delta}\}$.

Let $\overline{\tau}\in [0,1/2)$ be such that
\begin{equation}
\label{lm:avgr2}\Past\left(\lim_{t\rightarrow\infty}(t+1)^{\overline{\tau}}\left\|\mathbf{z}_{t}^{\delta}\right\|=0\right)=1
\end{equation}
for all $n$. Such a $\overline{\tau}$ always exists by Lemma~\ref{lm:avgc}. We now show that there exists $\tau$ such that $\overline{\tau}<\tau<1/2$ for which the claim holds. To this end, choose $\widetilde{\tau}\in (\tau,1/2)$ and let $\mu=1/2 (\overline{\tau}+\widetilde{\tau})$. For each $\delta>0$ recall the construction in~\eqref{lm:avgc2}-\eqref{lm:avgc8} and the $\mathcal{F}_{t}$-adapted process $\{\mathbf{z}_{t}^{\delta}\}$ satisfies
\begin{eqnarray}
\label{lm:avgr3}
\left\|\mathbf{z}_{t+1}^{\delta}\right\|^{2}& \leq & \left\|I_{M}-\alpha_{t}\Gamma_{t}^{\delta}\right\|^{2}\left\|\mathbf{z}_{t}^{\delta}\right\|^{2}+\alpha_{t}^{2}\left\|U_{t}^{\delta}\right\|^{2}+\alpha_{t}^{2}\left\|J_{t}^{\delta}\right\|^{2}+ 2\alpha_{t}\left(\mathbf{z}_{t}^{\delta}\right)^{T}\left(I_{M}-\alpha_{t}\Gamma_{t}^{\delta}\right)J_{t}^{\delta}\nonumber \\ & & +2\alpha_{t}\left\|U_{t}^{\delta}\right\|\left(\left\|I_{M}-\alpha_{t}\Gamma_{t}^{\delta}\right\|\left\|\mathbf{z}_{t}^{\delta}\right\|+\alpha_{t}\left\|J_{t}\right\|\right).
\end{eqnarray}
Since $\tau_{1}>\tau_{2}+1/(2+\Vap_{1})+1/2$, by Lemma~\ref{lm:estc} and~\eqref{lm:avgc3} the process $\{U_{t}^{\delta}\}$ may be chosen such that\footnote{For $\mathbb{R}_{+}$ valued sequences $\{f_{t}\}$ and $\{g_{t}\}$ the notation $f_{t}=o(g_{t})$ means that $f_{t}/g_{t}\rightarrow 0$ as $t\rightarrow\infty$. For stochastic sequences the $o(\cdot)$ is to be interpreted a.s. or pathwise.}
\begin{equation}
\label{lm:avgr4}
\left\|U_{t}^{\delta}\right\|=o\left((t+1)^{-1/2}\right).
\end{equation}
Since $\left\|\mathbf{z}_{t}^{\delta}\right\|=o\left((t+1)^{-\overline{\tau}}\right)$ (by hypothesis), we obtain
\begin{equation}
\label{lm:avgr400}
2\alpha_{t}\left\|U_{t}^{\delta}\right\|\left\|I_{M}-\alpha_{t}\Gamma_{t}^{\delta}\right\|\left\|\mathbf{z}_{t}^{\delta}\right\|=o\left((t+1)^{-3/2-\overline{\tau}}\right).
\end{equation}
The existence of the second moment of the observation noise process and the boundedness of $\{\mathcal{K}_{t}^{\delta}\}$ imply
\begin{equation}
\label{lm:avgr5}
\Past\left(\lim_{t\rightarrow\infty}(t+1)^{-1/2-\Vap}\left\|J_{t}^{\delta}\right\|=0\right)=1
\end{equation}
for each $\Vap>0$ and, hence
\begin{equation}
\label{lm:avgr6}
2\alpha_{t}^{2}\left\|U_{t}^{\delta}\right\|\left\|J_{t}^{\delta}\right\|=o\left((t+1)^{-3/2-\overline{\tau}}\right).
\end{equation}
Since $2\mu=\overline{\tau}+\widetilde{\tau}$ and $\widetilde{\tau}<1/2$, by~\eqref{lm:avgr5} we note that
\begin{equation}
\label{lm:avgr7}
\sum_{t\geq 0}(t+1)^{2\mu}\alpha_{t}\left\|U_{t}^{\delta}\right\|\left\|I_{M}-\alpha_{t}\Gamma_{t}^{\delta}\right\|\left\|\mathbf{z}_{t}^{\delta}\right\|<\infty.
\end{equation}
Similarly we have
\begin{equation}
\label{lm:avgr8}
\sum_{t\geq 0}(t+1)^{2\mu}\alpha_{t}^{2}\left\|U_{t}^{\delta}\right\|\left\|J_{t}^{\delta}\right\|<\infty,~~~\sum_{t\geq 0}(t+1)^{2\mu}\alpha_{t}^{2}\left\|U_{t}^{\delta}\right\|^{2}<\infty.
\end{equation}
Now consider the terms $\alpha_{t}^{2}\|J_{t}^{\delta}\|^{2}$. Since the second moment of the observation noise process exists, $\{\mathcal{K}_{t}^{\delta}\}$ is uniformly bounded and $2\mu<1$, it can be shown that
\begin{equation}
\label{lm:avgr9}
\sum_{t\geq 0}(t+1)^{2\mu}\alpha_{t}^{2}\|J_{t}^{\delta}\|^{2}<\infty.
\end{equation}
Now let $\{W_{t}^{\delta}\}$ denote the $\mathcal{F}_{t+1}$ sequence given by
\begin{equation}
\label{lm:avgr10}
W_{t}^{\delta}=\alpha_{t}\left(\mathbf{z}_{t}^{\delta}\right)^{T}\left(I_{M}-\alpha_{t}\Gamma_{t}^{\delta}\right)J_{t}^{\delta}.
\end{equation}
We note that $\East[W_{t}^{\delta}|\mathcal{F}_{t}]=0$ for all $t$ and (at least for $t$ large) we have $\East[\left(W_{t}^{\delta}\right)^{2}|\mathcal{F}_{t}]\leq\alpha_{t}^{2}\left\|\mathbf{z}_{t}^{\delta}\right\|^{2}\left\|J_{t}^{\delta}\right\|^{2}$.
Since the second moment of the observation noise process exists and $\{\mathcal{K}_{t}^{\delta}\}$ is uniformly bounded, we obtain
\begin{equation}
\label{lm:avgr12}
\East\left[\left(W_{t}^{\delta}\right)^{2}~|~\mathcal{F}_{t}\right]=o\left((t+1)^{-2-2\overline{\tau}}\right).
\end{equation}
Hence
\begin{equation}
\label{lm:avgr13}
\East\left[(t+1)^{4\mu}\left(W_{t}^{\delta}\right)^{2}~|~\mathcal{F}_{t}\right]=o\left((t+1)^{-2-2\overline{\tau}+4\mu}\right)=o\left((t+1)^{-2+2\widetilde{\tau}}\right).
\end{equation}
Since $2\widetilde{\tau}<1$, the sequence on the left hand side of~\eqref{lm:avgr13} is summable and by Proposition~\ref{prop:DF} we conclude that $\sum_{t\geq 0}(t+1)^{2\mu}W_{t}^{\delta}$ exists and is finite. Since $\Gamma_{t}^{\delta}\rightarrow I_{M}$ uniformly and $\alpha_{t}\rightarrow 0$ as $t\rightarrow\infty$, we have
\begin{equation}
\label{lm:avgr15}
\left\|I_{M}-\alpha_{t}\Gamma_{t}^{\delta}\right\|^{2}\leq \left(1-a(t+1)^{-1}\right)
\end{equation}
for all $t$ large enough. Thus (eventually) we have from~\eqref{lm:avgr3}
\begin{equation}
\label{lm:avgr16}
\left\|\mathbf{z}_{t+1}^{\delta}\right\|^{2}\leq\left(1-a(t+1)^{-1}\right)\left\|\mathbf{z}_{t}^{\delta}\right\|^{2}+d_{t}(t+1)^{-2\mu}
\end{equation}
where the term $d_{t}(t+1)^{-2\mu}$ corresponds to all the residuals. Moreover by~\eqref{lm:avgr4}-\eqref{lm:avgr15} the limit $\lim_{t\rightarrow\infty}\sum_{s=0}^{t}d_{s}$ exists and is finite. Since $a\geq 1>2\mu$, an immediate application of Proposition~\ref{prop:Fab-1} yields
\begin{equation}
\label{lm:avgr18}
\limsup_{t\rightarrow\infty}(t+1)^{2\mu}\left\|\mathbf{z}_{t}^{\delta}\right\|^{2}<\infty~~\mbox{a.s.}
\end{equation}
Hence, there exists $\tau$ with $\overline{\tau}<\tau<\mu$, such that
$(t+1)^{\tau}\left\|\mathbf{z}_{t}^{\delta}\right\|\rightarrow 0$ a.s. as $t\rightarrow\infty$. Since the above holds for all $\delta>0$, we conclude that
$(t+1)^{\tau}\left\|\mathbf{z}_{t}\right\|\rightarrow 0$ a.s. as $t\rightarrow\infty$. Thus, for every $\overline{\tau}$ for which the convergence in~\eqref{lm:avgr1} holds there exists $\tau\in (\overline{\tau},1/2)$ for which the convergence continues to hold. Hence, by induction we conclude that the required convergence holds for all $\tau\in [0,1/2)$.
\end{proof}
}

\bibliographystyle{unsrt}
\bibliography{CentralBib}

\end{document}